\pdfoutput=1
\RequirePackage{ifpdf}
\ifpdf % We~are running pdfTeX in pdf mode
\documentclass[pdftex]{sigma}
\else
\documentclass{sigma}
\fi

\usepackage{bm}

\usepackage[subsetnonamb, bbsets]{jkmath}

\usepackage{breakurl}

\newtheorem{Theorem}{Theorem}[section]

\newtheorem{Lemma}[Theorem]{Lemma}
\newtheorem{Proposition}[Theorem]{Proposition}
\newtheorem{Corollary}[Theorem]{Corollary}

\theoremstyle{definition}
\newtheorem{Definition}[Theorem]{Definition}
\newtheorem{Example}[Theorem]{Example}
\newtheorem{Remark}[Theorem]{Remark}

\DeclareMathOperator{\ct}{ct}

\numberwithin{equation}{section}

\begin{document}
\allowdisplaybreaks

\newcommand{\arXivNumber}{2310.17362}

\renewcommand{\PaperNumber}{057}

\FirstPageHeading

\ShortArticleName{Intermediate Macdonald Polynomials and Their Vector Versions}

\ArticleName{Intermediate Macdonald Polynomials\\ and Their Vector Versions}

\Author{Philip SCHL\"OSSER}

\AuthorNameForHeading{P.~Schl\"osser}

\Address{IMAPP-Mathematics, Radboud University, Heyendaalseweg 135, \\ 6525 AJ Nijmegen, The Netherlands}
\Email{\mail{philip.schloesser@ru.nl}}
\URLaddress{\url{https://philip98.github.io}}

\ArticleDates{Received October 30, 2025, in final form May 13, 2026; Published online June 04, 2026}

\Abstract{Intermediate Macdonald polynomials for an affine root system $S$ with fixed origin and finite Weyl group $W_0$ are orthogonal polynomials invariant under a parabolic subgroup $W_J\le W_0$. The extreme cases of $W_J=1$ and $W_J=W_0$ correspond to the non-symmetric and symmetric Macdonald polynomials, respectively. In this paper, we use double-affine Hecke algebras to study their basic properties, including that they form an orthogonal basis and that they diagonalise a commutative algebra of difference-reflection operators, and calculate their norms. Finally, we provide two interpretations of intermediate Macdonald polynomials as vector-valued polynomials and connect them to the literature.}


\Classification{33D52; 33D80; 20C08}

\section{Introduction}
Let $R$ be a (finite) root system and $L$ a lattice associated with $R$, such
as the weight lattice, the coweight lattice, or the coroot lattice. Let
$A:= K[L]$ be the group algebra of $L$ over a suitable field $K$. Orthogonal
polynomials associated with $R$, such as the ones detailed in~\cite{cherednikDAHA,macdonald,opdam}, are typically elements of $A$ that are either
symmetric under the Weyl group $W_0$ of $R$ (or transform according to a
given 1-dimensional representation) or possess no symmetry at all.

In this paper, we study the intermediate cases: polynomials
that transform under a character~$\epsilon$ of a parabolic subgroup $W_J\le W_0$, in the
context of Macdonald polynomials.

In the context of Heckman--Opdam polynomials, the $W_J$-invariant polynomials
have been considered in~\cite{mvpVectorValued}. Our approach is
different and leads
to the notion of \emph{intermediate Macdonald polynomials}.

In the case where $R$ is of type $A_n$ and $\epsilon=1$,
these polynomials have been studied extensively, e.g., in
\cite{BDF97, baratta,goodberry,lapointe} under the name ``Macdonald polynomials with prescribed symmetry'', and connected to the
equivariant K-theory of parabolic flag Hilbert schemes in~\cite{goodberryGeometric} as well as to the theory of parahoric Whittaker functions in~\cite{BBBG24}. In a preprint~\cite{FKhM} that appeared shortly after this one, it is furthermore shown that certain specialisations of intermediate Macdonald polynomials can be obtained as characters of representations of parahoric subalgebras of affine Kac--Moody algebras.

For general Macdonald data and characters $\epsilon$ of $W_J$, we prove basic facts about the intermediate Macdonald polynomials,
such as leading terms, orthogonality, and that they diagonalise a~commutative algebra
of difference-reflection operators, and compute their norms.

Finally, we provide two vector-valued interpretations for $\epsilon=1$,
and connect them to the literature.
For the first interpretation as vector-valued polynomials, we can artificially increase the symmetry of elements of $A_J:= A^{W_J}$
by tensoring with the finite-dimensional vector space~${K[W_0/W_J]}$. In fact,
\[
 A_J \cong (A\otimes K[W_0/W_J])^{W_0}
\]
as modules over $A_0:= A^{W_0}$ by mapping $A_J\ni f\mapsto f\otimes e(W_J)$ and then symmetrising ($e(W_J)$ is the basis vector of $K[W_0/W_J]$ corresponding to the coset of the neutral element). However, in the context of Macdonald
polynomials, the group $W_0$ and its subgroups are rather artificial,
and it makes much more sense to phrase our constructions in terms of
the double-affine Hecke algebra $\tilde{\mathfrak{H}}$ from
\cite{macdonald}. Luckily, the authors of~\cite{stokmanInduced} show that it is possible to establish a $\tilde{\mathfrak{H}}$-module structure (even more than one) on $A\otimes K[W_0/W_J]$, which they call
$\mathbb{M}_J(\phi)$. We show that the ``spherical'' vectors
of $\mathbb{M}_J(\phi)$ are indeed in $A_0$-linear equivalence with $A_J$.
This result has been proven independently in~\cite{Ven25}.

This vector-valued interpretation is especially valuable if we seek
to connect our theory to representation theory. In the classical theory
(i.e., Lie groups),
$(A\otimes K[W_0/W_J])^{W_0}$ is exactly the form that suitable
matrix-spherical functions take. In fact, besides the well-known
correspondence between symmetric Heckman--Opdam polynomials ($W_J=W_0$)
and zonal spherical functions of symmetric pairs from, e.g., \cite{heckmanSchlichtkrull}, a concrete connection between intermediate
Heckman--Opdam polynomials and matrix spherical functions of symmetric
pairs has been shown in~\cite{mvpBC1,mvpVectorValued} for some examples.
For quantum groups, a similar connection between symmetric Macdonald polynomials
and zonal spherical functions of quantum symmetric pairs was proven
in~\cite{letzter}, so it stands to reason to look for matrix spherical
functions in the Macdonald version of~${(A\otimes K[W_0/W_J])^{W_0}}$, i.e., the spherical elements of~${\mathbb{M}_J(\phi)}$.

For a second interpretation of $A_J$ as a ring of vector-valued
polynomials, recall that by a~well-known result by Chevalley,
$A_0$ is a free (commutative) $K$-algebra,
i.e., a polynomial ring. The same, however, cannot be said about $A_J$
(except when $W_J=W_0$) or even $A$. Thus, we are confronted with the
pedantic philosophical question of how exactly intermediate and even non-symmetric
Macdonald polynomials can even be thought of as polynomials in
the first place. A possible solution to this conundrum comes from noting
that $A_J$ is a free $A_0$-module (as is detailed in~\cite{pittieSteinberg}), so that we
can view elements of $A_J$ as column vectors of elements of $A_0$, i.e.,
as vector-valued polynomials. The inner product on $A_J$ can then be
rewritten in terms of a matrix weight with symmetric entries, which
allows for an interpretation in the theory of vector-valued orthogonal
polynomials.

In Section~\ref{sec:daha-general}, we begin by introducing some notions about double-affine Hecke algebras for arbitrary
affine root systems $S$, as well as setting up the necessary notions of duality, labels, and
the monomial order on the ring where our orthogonal polynomials live. The notation (with
exception of the $\tau$'s, as explained in Remark~\ref{sec:rem-taus}\,(ii)) and
ordering of topics have been taken from~\cite{macdonald}, which we will refer to for most of the proofs.

In Section~\ref{sec:macdo-parabolic}, we introduce various tools that are specifically related to parabolic subgroups,
and that are not taken wholesale from~\cite{macdonald}. They include an adaption of the symmetrisers from~\cite[Section~5.5]{macdonald} as well as a study of $W_J$-orbits of monomials.

Section~\ref{sec:poly} is the central section of this paper and contains its protagonist: the ($\epsilon$-symmetric) intermediate
Macdonald polynomials. We introduce them in analogy with the $\epsilon$-symmetric Macdonald polynomials
from \cite[Section~5.7]{macdonald} by applying a suitable symmetriser to non-symmetric Macdonald polynomials
and then normalising. We show that the symmetric versions form a~basis of the $W_J$-symmetric polynomials,
and that they are orthogonal. Then we proceed to relate their norms to those of the non-symmetric
Macdonald polynomials.

In Section~\ref{sec:jasper}, we consider the first
vector-valued interpretation of $A_J$. We introduce the $Y$-parabolically
induced modules $\mathbb{M}_J(\phi)$ from \cite[Section~3.2]{stokmanInduced} and then prove that their spaces of spherical vectors
is in $A_0$-linear equivalence with $A_J$.

Finally, in Section~\ref{sec:w0-invariant}, we consider the second vector-valued interpretation
and show the existence of a symmetric matrix weight.
To illustrate this general theory, we then consider two examples: the case
of non-symmetric Askey--Wilson polynomials
\big(type $\bigl(C_1^\vee, C_1\bigr)$ with $W_J=1$, so that $A^{W_J}=A$\big) from~\cite{koornwinderBouzeffour}, and the case
of type $A_2$ with a nontrivial parabolic subgroup~${W_J\le W_0}$ that was also considered in the classical case
in~\cite{mvpVectorValued}.

\section{Double-affine Hecke algebras}\label{sec:daha-general}
For the convenience of the reader, we now recap the theory from~\cite{macdonald} that we will rely on in
the following sections. Readers wishing to skip the recap are advised to consult Remark~\ref{sec:rem-taus}
for notational differences to Macdonald, and then proceed to Section~\ref{sec:macdo-parabolic}.

Let $E$ be a finite-dimensional Euclidean space with translation space $V$, and let $F$ be the
real vector space of affine functions $E\to\R$. Let $S\subset F$ be an affine root system (see~\cite{macdoARS}) and~$W_S$ its affine Weyl group.

\subsection{Duality and labels}
To each possible choice of $S$ (potentially up to embedding in a larger
root system, cf.\ the remark immediately proceeding \cite[Section~1.4]{macdonald}), we are going to associate another affine root system~$S'$, two
(linear) root systems $R$, $R'$, and two lattices $L$, $L'$.

\begin{Lemma}\label{sec:lem-all-affine-root-systems}
 Let $S$ be irreducible, then after choosing an appropriate origin to identify $E=V$ and after appropriate re-scaling,
 one of the following is true:
 \begin{enumerate}\itemsep=0pt
 \item[$(i)$] There is a reduced linear root system $R\subset V$ such that $S=S(R)$.
 \item[$(ii)$] There is a reduced linear root system $R\subset V$ such that $S=S(R)^\vee$.
 \item[$(iii)$] There is an affine root system $\tilde{S}$ of type $\bigl(C_n^\vee, C_n\bigr)$ of which $S$ forms
 a $W_{\tilde{S}}$-invariant subsystem.
 \end{enumerate}
 Here, $S(R)$ is as in {\rm \cite[Section~1.2.1]{macdonald}}.
\end{Lemma}
\begin{proof}
 This is the content of \cite[Section~1.3]{macdonald}.
 For reference, we mention for (iii), which orbits of $\tilde{S}$ make up
 $S$.
 We adopt the notation of \cite[Section~1.3]{macdonald}
 and define $O_1,\dots,O_5$ to be the orbits satisfying
$
 a_n\in O_1$, $ 2a_n\in O_2$, $ a_0\in O_3$, $ 2a_0\in O_4$, $ a_1,\dots,a_{n-1}\in O_5$
 ($a_0,\dots,a_n$ as in \cite[Section~1.3.18]{macdonald}).
 Then $S$ consists of the following
 $\tilde{S}$-orbits:
 \begin{gather*}
 BC_n\colon\ O_1,\ O_4,\ O_5,\qquad
 (BC_n, C_n)\colon\ O_1,\ O_2,\ O_4,\ O_5,\qquad
 \bigl(C_n^\vee, BC_n\bigr)\colon\ O_1,\ O_2,\ O_3,\ O_5,\\
 \bigl(B_n, B_n^\vee\bigr)\colon\ O_1,\ O_2,\ O_5,\qquad
 \bigl(C_n^\vee, C_n\bigr)\colon\ O_1,\ O_2,\ O_3,\ O_4,\ O_5.\tag*{\qed}
 \end{gather*} \renewcommand{\qed}{}
\end{proof}

For the remainder of this paper, we replace $S$ by $\tilde{S}$ if Lemma~\ref{sec:lem-all-affine-root-systems}\,(c) is satisfied.

\begin{Definition}
 If $S$ is irreducible, and
 \begin{enumerate}\itemsep=0pt
 \item[(i)] $S=S(R)$ for a reduced (linear) root system $R$, let
$
 S':= S\bigl(R^\vee\bigr)$, $ R':= R^\vee$, $ L:= P(R)$, $ L':= P\bigl(R^\vee\bigr)$.
 Furthermore, for $\alpha\in R$ write
 $\alpha':=\alpha^\vee\in R'$.
 \item[(ii)] $S=S(R)^\vee$ for a reduced (linear) root system $R$, let
$
 S':= S$, $ R':= R$, $ L:=L':=P\bigl(R^\vee\bigr)$.
 Furthermore, for $\alpha\in R$ write
 $\alpha':=\alpha\in R'$.
 \item[(iii)] $S$ is of type $\bigl(C_n^\vee, C_n\bigr)$, let $R$ be the linear root system of
 type $C_n$ such that $S=S(R)\cup S(R)^\vee$. Then let
$
 S':= S$, $ R':= R$, $ L:=L':= \Z R^\vee$.
 Furthermore, for $\alpha\in R$ write
 $\alpha':=\alpha\in R'$.
 \end{enumerate}
 If $S=S_1\cup\cdots\cup S_n$ is reducible, assume without loss of generality that
 every $S_i$ is either $S(R_i)$ ($R_i$ reduced), $S(R_i)^\vee$ ($R_i$ reduced), or
 that $S(R_i)\cup S(R_i)^\vee$ ($R_i$ of type $C_{n_i}$). This has the side-effect that all
 constant affine functions contained in $S_i$ and $S'_i$, are half-integer-valued. Write
 $c$ for the constant function that is $1$ everywhere.
 Then define $S'_i$, $R_i$, $R'_i$, $L_i$, $L'_i$ as above and take
 \begin{gather*}
 S':= S'_1\cup\cdots\cup S'_n,\qquad
 R:= R_1\cup\cdots\cup R_n,\qquad
 R':= R'_1\cup\cdots\cup R'_n,\\
 L:= L_1\oplus\cdots\oplus L_n,\qquad
 L':= L'_1\oplus\cdots\oplus L'_n.
 \end{gather*}
 We interpret $R$, $R'$, $L$, $L'$ as subsets of $V$ and $S$, $S'$ as sets of
 affine-linear functions on $V$, and we call
 $(S,S',R,R',L,L')$ a \emph{set of duality data}.
\end{Definition}

\begin{Remark}\qquad
\begin{enumerate}\itemsep=0pt
 \item[(i)]
 If $(S,S',R,R',L,L')$ is a set of duality data, then so is
 $(S',S,R',R,L',L)$.
 \item[(ii)]
 Let $\langle\cdot,\cdot\rangle$ be the inner product of the translation
 space $V$, and use it to identify $V\cong V^*$. If $D\colon F\to V$ is the gradient map, i.e.,
 $f(p+v)=f(p) + \langle Df, v\rangle$ for $f\in F$, $p\in E$, $v\in V$, we then have $DS\subset L$, $DS'\subset L'$. Furthermore, we have
 have $R^\vee\subset S'$ and $R^{\prime\vee}\subset S$.
 \item[(iii)]
 Furthermore, $\langle R, L'\rangle,\langle R',L\rangle\subset\Z$. In
 particular, unless $S$ contains components of type $\bigl(C_n^\vee, C_n\bigr)$, the pairings of $R\times L'$, $R'\times L$ are
 perfect pairings. If $S$ does contain a component of type $\bigl(C_n^\vee, C_n\bigr)$ and if $a\in S$ lies in the orbit of the
 highest root of $R$, we have $\langle a, L'\rangle=\langle a',L\rangle=2\Z$.
 \item[(iv)]
 It is possible to choose systems of simple roots $(a_i)_{i\in I}$, $(a'_i)_{i\in I}$, and
 $(\alpha_i)_{i\in I_0}$ of $S$, $S'$ and $R$, respectively, such that $\alpha_i^\vee=a_i'$ ($i\in I_0$). Here, $I_0\subset I$ indexes the
 simple linear roots and leaves out only one (non-linear) affine root per
 irreducible component. For each irreducible component $R_i$, depending on whether $R'_i=R_i$ or $R'_i=R_i^\vee$, write
 $\alpha':=\alpha$ or $\alpha':=\alpha^\vee$ for~${\alpha\in R_i}$.
\end{enumerate}
\end{Remark}

\begin{Definition}
 Let $(S,S',R,R',L,L')$ be a set of duality data, write
 $W_0$ for the (finite) Weyl group of $R$ and $R'$ (they have the same).
 Then define the \emph{extended affine Weyl groups}
 \begin{align*}
 W:= W(R,L') := W_R\ltimes t(L'),\qquad
 W':= W(R',L).
 \end{align*}
 Here, $t(\lambda')\colon V\to V$ is the translation $x\mapsto x + \lambda'$
 (for $\lambda'$). Since $W$ acts on $V$ by affine transformations,
 it also acts on $F$, the space of affine functions. In particular,
 $W$ acts on~$L'$ and $S$; similarly, $W'$ acts on $L$ and $S'$. Furthermore, we have
 $W_S\le W$ and $W_{S'}\le W'$, where~$W_S$ (resp.\ $W_{S'}$)
 is the reflection group generated by $s_i=s_{a_i}$ (resp.\ $s_{a'_i}$)
 for $i\in I$ ($s_a$ is the affine reflection in $a$, see \cite[Section~1.1.6]{macdonald}).

 Write $W_0:= W_{R}=W_{R'}$ for the reflection group generated by $s_i = s_{\alpha_i}$ for $i\in I_0$. It acts on $S$, $S'$, $R$, $R'$, $L$, $L'$ and is
 the parabolic subgroup of $W_S$ and $W_{S'}$ generated by the
 index set $I_0$.
\end{Definition}

\begin{Definition}
 A function $k\colon S\to \R$ is called a \emph{$W$-labelling} if $k$ is constant on
 $W$-orbits.
\end{Definition}

\begin{Remark}
Let $S$ be irreducible.
\begin{enumerate}\itemsep=0pt
 \item[(i)]
 If $S=S(R)$ for $R$ simply-laced, there is only one $W$-orbit, and hence any
 labelling $k$ is given by one number.
 \item[(ii)]
 If $S=S(R)$ or $S(R)^\vee$ for $R$ not simply-laced, there are two $W$-orbits,
 so any labelling $k$ is given by two numbers $k_s$, $k_l$ for short and long roots,
 respectively.
 \item[(iii)]
 If $S$ is of type $(C_n,C_n^\vee)$, we have $W=W_S$ and therefore we have
 the four $(n=1)$ or five orbits $O_1,\dots,O_5$ from the proof of Lemma~\ref{sec:lem-all-affine-root-systems}. In that case,
 a labelling $k$ is given by five numbers $k_1,\dots,k_5$.
\end{enumerate}
\end{Remark}

\begin{Definition}
 The \emph{dual labelling} of a $W$-labelling $k$ is a $W'$-labelling $k'\colon S'\to\R$
 given by
 \begin{enumerate}\itemsep=0pt
 \item[(i)] $S=S(R)\colon k'(\alpha^\vee+c):= k(\alpha+rc)$ for $\alpha\in R, r\in\frac{1}{2}\Z$;
 \item[(ii)] $S=S(R)^\vee\colon k':= k$;
 \item[(iii)] Type $\bigl(C_n^\vee, C_n\bigr)\colon k'_5:= k_5$ and
 \[
 \mqty(k'_1\\k'_2\\k'_3\\k'_4) :=
 \frac{1}{2}\mqty(1 & 1 & 1 & 1\\1 & 1 & -1 & -1\\1 & -1 & 1 & -1\\1 & -1 & -1 & 1)\mqty(k_1\\k_2\\k_3\\k_4)
 \]
 \end{enumerate}
 on irreducible components.
\end{Definition}
Evidently, if $k'$ is dual to $k$, then $k$ is dual to $k'$.

\begin{Lemma}\label{sec:lem-ks1}
 Let $S_k\subset S$ be an irreducible component and let $i\in I_0$ such that $a_i\in S_k$.
 Then
 \[
 k(a_i) + k(2a_i) = k'(a'_i) + k'(2a'_i).
 \]
\end{Lemma}
\begin{proof}
 We differentiate between the three cases:
 \begin{enumerate}\itemsep=0pt
 \item If $S_k=S(R_k)$, then $a_i=\alpha_i$ and $a'_i=\alpha_i^\vee$,
 so that
$
 k'(a'_i)=k'\bigl(\alpha_i^\vee\bigr)=k(\alpha_i)=k(a_i)
$
 and~${k'(2a'_i)=0=k(2a_i)}$.
 \item If $S_k=S(R_k)^\vee$, then $a_i'=a_i$ and $k=k'$, so that
 $k'(a'_i)=k(a_i)$ and $k'(2a'_i)=0=k(2a_i)$.
 \item If $S_k$ is of type $\bigl(C_n^\vee, C_n\bigr)$,
 then $a_i$ is indivisible and linear. Hence $a_i$ belongs either to
 $O_1$ or to $O_5$. In the first case we have $2a_i\in O_2$ and hence
 \[
 k'(a'_i) + k'(2a'_i)
 = k'_1 + k'_2 = k_1 + k_2
 = k(a_i) + k(2a_i).
 \]
 In the latter case, $2a_i\not\in S$ and hence
 \[
 k'(a'_i) + k'(2a'_i) = k'_5 = k_5 = k(a_i) + k(2a_i).\tag*{\qed}
 \]
 \end{enumerate}\renewcommand{\qed}{}
\end{proof}

\begin{Lemma}\label{sec:lem-ks2}
 Let $S_k\subset S$ be an irreducible component and let $i\in I\setminus I_0$ be such that $a_i\in S_k$ $($in particular, $a_i$ is not linear$)$.
 Let $j\in I_0$ be such that $W_0\alpha_j$ contains the highest root
 $\phi$ of $R_k$. Then
$
 k(a_i)+k(2a_i) = k'(a'_j) - k'(2a'_j)$.
\end{Lemma}
\begin{proof}
 We distinguish between the three cases.
 \begin{enumerate}\itemsep=0pt
 \item If $S_k=S(R_k)$, then $Da_i=-\phi$, so that
$
 k(a_i)=k(-\phi) = k(\phi) = k'\bigl(\phi^\vee\bigr)$.
 Since $\alpha_j\in W_0\phi$, we have $\alpha_j^\vee=a'_j\in W_0\phi^\vee$ and hence $k(a_i)=k'(a'_j)$. Since $S_k$ is reduced,
 the elements $2a_i$, $2a'_j$ are not in $S$.
 \item If $S_k=S(R_k)^\vee$, then $Da_i=-\phi^\vee$, so that
 \[
 k(a_i) = k\bigl(-\phi^\vee\bigr)=k\bigl(\phi^\vee\bigr)=
 k'\bigl(\phi^\vee\bigr).
 \]
 Since $\alpha_j\in W_0\phi$, we again have $a'_j\in W_0\phi^\vee$, so that $k(a_i)=k'(a'_j)$.
 \item If $S_k$ is of type $\bigl(C_n^\vee, C_n\bigr)$,
 we have $a_i\in O_3$, so that
 \[
 k(a_i)+k(2a_i)=k_3+k_4 = k'_1-k'_2.
 \]
 Since $\alpha_j$ is long, its corresponding simple $a'_j$ (of $S_k$) is short (see \cite[Section~1.4.3]{macdonald}, where~${j=n}$)
 and lies in $O_1$.
 Therefore,
 \[
 k(a_i)+k(2a_i) = k'_1-k'_2 = k'(a'_j)-k'(2a'_j).\tag*{\qed}
 \]
 \end{enumerate}\renewcommand{\qed}{}
\end{proof}

\subsection{Extended affine Weyl group}
Write
$
 S_1:=\set{a\in S\where 2a\not\in S}
$
for the subsystem of inmultipliable roots. We shall now have a~closer look
at the extended affine Weyl groups $W$, $W'$, in particular at their length functions.

\begin{Definition}
 For $w\in W$, define
 \[
 S_1(w):=S^+\cap w^{-1} S^-
 = \set{a\in S^+\where wa\in S^-},
 \]
 which we can also interpret as the set of hyperplanes lying between the fundamental
 alcove $C$ and $w^{-1}C$. Define the \emph{length} of $W$ to be
$
 \ell(w) := \# S_1(w)$.
 Write $\Omega:=\set{w\in W\where \ell(w)=0}$.
\end{Definition}

\begin{Lemma}[{\cite[below equation~(2.2.3)]{macdonald}}]
 The group $\Omega$ is a finite abelian group, and we have
 $W=W_S\rtimes\Omega$.
\end{Lemma}

\begin{Lemma}[{\cite[Section~2.2.4]{macdonald}}]\label{sec:lem-lengths-add}
 If $w,v\in W$, then $\ell(vw)\le \ell(v)+\ell(w)$ with equality iff~${S_1(w)\subset S_1(vw)}$.
\end{Lemma}

\begin{Lemma}[{\cite[Section~2.2.8]{macdonald}}]\label{sec:lem-length-change}
 Define $\sigma\colon S\to\R$ by $\sigma(a):=\pm1$ for $a\in S^\pm$. For all
 $w\in W$ and $i\in I$, we have
 \begin{gather*}
 \ell(s_iw) = \ell(w) + \sigma\bigl(w^{-1}a_i\bigr),\qquad
 \ell(ws_i) = \ell(w) + \sigma(wa_i).
 \end{gather*}
\end{Lemma}

\begin{Lemma}[{\cite[Section~2.2.9]{macdonald}}]\label{sec:lem-flipped-roots}
 Every element $w\in W$ can be expanded as a reduced expression of the form
 $w=us_{i_1}\cdots s_{i_p}$, where $p=\ell(w)$ and $u\in\Omega$. For such a $w$, we have
$ S_1(w) = \set{b_1,\dots,b_p}$,
 where $b_r := s_{i_p}\cdots s_{i_{r+1}} (a_{i_r})$.
\end{Lemma}

\begin{Definition}
 $W_S$ as a Coxeter group comes equipped with a Bruhat order. We extend it to~$W$ as follows: $uw\le vw'$ for $u,v\in\Omega$ and $w,w'\in W_S$ iff
 $u=v$ and $w\le w'$.
\end{Definition}

\begin{Lemma}\label{sec:lem-existence-shortest-coset-dec}
 Let $W_J\le W_S$ be a parabolic subgroup. Define
 \[
 W^J := \set{v\in W\where\forall w\in W_J\colon \ell(vw)\ge\ell(v)},
 \]
 the set of shortest coset representatives. Then every element
 $w\in W$ has a unique decomposition~$vw'$ with $v\in W^J$, $w'\in W_J$ and $\ell(w)=\ell(v)+\ell(w')$.
\end{Lemma}
\begin{proof}
 If $W$ is a Coxeter group, the claim is a well-known result that
 can be found, e.g., in \cite[Section~5.12, Proposition~1.12]{humphreysRefl}. Since we are here also allowing for the semidirect
 product with the finite automorphism group $\Omega$, there is something left to prove.

 Write $w=u\tilde{w}$ for $u\in\Omega,\tilde{w}\in W_S$. By \cite[Section~5.12, Proposition~1.12]{humphreysRefl}, $\tilde{w}=\tilde{v}\tilde{w}'$
 for~${\tilde{w}'\in W_J}$ and $\tilde{v}$ a shortest coset representative. Then
 picking $w':= \tilde{w}'$ and $v:=u\tilde{v}$ does the trick: evidently,
 $w'=\tilde{w}'\in W_J$, and the map $W_S\to W, x\mapsto ux$ is an order-isomorphism
 onto $uW_S$, both with respect to the quasi-order induced by the length function
 and the Bruhat order.
\end{proof}

\subsection[rho, u, v, and r]{$\boldsymbol{\rho}$, $\boldsymbol{u}$, $\boldsymbol{v}$, and $\boldsymbol{r}$}
We shall now decompose translations in terms of shortest coset
representatives with respect to~$W_0$ and define a few auxiliary objects
that make use of the $W$-labelling $k$.

\begin{Definition}
 Define
 \[
 \rho_{k'} := \frac{1}{2}\sum_{\alpha\in R} k'\bigl(\alpha^\vee\bigr)\alpha
 \qquad \rho'_k := \frac{1}{2}\sum_{\alpha'\in R'} k\bigl(\alpha^{\prime\vee}\bigr)\alpha' \in V.
 \]
\end{Definition}

\begin{Lemma}\label{sec:lem-rho-inner-prod}
 We have $\big\langle\rho_{k'},\alpha_i^\vee\big\rangle=k'\bigl(\alpha_i^\vee\bigr)$
 for $i\in I_0$. In particular,
 if $k'(a')\ne0$ for all~${a'\in S'}$, then $\rho_{k'}$ is regular under
 the action of $W_0$ $($viewed as an element of $E\cong V)$. If~${k'(a')>0}$ for all
 $a'\in S'$, it is even strictly dominant.
\end{Lemma}
\begin{proof}
 For $i\in I_0$, we have
 \[
 s_i \rho_{k'} = \frac{1}{2}\sum_{\alpha\in R^+} k'\bigl(\alpha^\vee\bigr)s_i\alpha
 = \frac{1}{2}\sum_{\alpha\in s_iR^+} k'\bigl(\alpha^\vee\bigr)\alpha.
 \]
 By \cite[Chapitre~VI, Section~1.6, Corollaire~1]{bourbaki}, the simple reflection $s_i$
 permutes the positive roots that are not $\alpha_i$. Therefore,
 we have $s_i R^+ = R^+\setminus\set{\alpha_i}\cup\set{-\alpha_i}$, hence
 \[
 s_i \rho_{k'} = \rho_{k'} - k'\bigl(\alpha_i^\vee\bigr) \alpha_i
 = \rho_{k'} - \big\langle \rho_{k'}, \alpha_i^\vee\big\rangle \alpha_i.
 \]
 Consequently, $k'\bigl(\alpha_i^\vee\bigr)=\langle\rho_{k'},\alpha_i^\vee\rangle$. So if
 $k'$ is never zero, $\rho_{k'}$ lies in the interior of a Weyl chamber. If $k'$
 is in particular positive, $\rho_{k'}$ lies in the interior of the fundamental
 Weyl chamber.
\end{proof}

\begin{Definition}
 For $\lambda'\in L'$, let $t(\lambda') =: u(\lambda')v(\lambda')\in W^0W_0$ be the decomposition from Lemma~\ref{sec:lem-existence-shortest-coset-dec}. Similarly, decompose $t(\lambda) =: u'(\lambda)v(\lambda)\in (W')^0W_0$
 for $\lambda\in L$.
\end{Definition}

 Define
 \[
 L'_+:=\set{\lambda'\in L'\where\forall i\in I_0\colon \langle\lambda',a_i\rangle\ge0}
 \]
 and analogously, $L_+$, and call the elements of $L'_+$, $L_+$ \emph{dominant}. The dominant elements form a fundamental domain for
 the action of $W_0$ on $L'$, $L$, respectively. Elements of $-L'_+$, $-L_+$ are called
 \emph{antidominant}. For $\lambda'\in L'$, write $\lambda'_+$, $\lambda'_-$
 for the unique dominant (resp.\ antidominant) element in $W_0\lambda'$.

\begin{Lemma}\label{sec:lem-v-antidominant}
 Let $\lambda'\in L'$ $($resp.\ $\lambda\in L)$, then $v(\lambda')\lambda'\in -L'_+$
$($resp.\ $v(\lambda)\lambda\in -L_+)$.
\end{Lemma}
\begin{proof}
 In \cite[Section~2.4]{macdonald}, this follows from the definition of $v(\lambda')$. That the definition used there is equivalent to ours can be seen
 from \cite[Section~2.4.9]{macdonald}.
\end{proof}

\begin{Corollary}
 We have $u(\lambda')=v(\lambda')^{-1}t(\lambda'_-)$ for $\lambda'\in L'$, where
 $\lambda'_-$ is the antidominant element in its $W_0$-orbit.
\end{Corollary}
\begin{proof}
 By Lemma~\ref{sec:lem-v-antidominant}, we have $v(\lambda')\lambda'=\lambda'_-$,
 so that
 \[
 v(\lambda')^{-1}t(\lambda'_-) = v(\lambda')^{-1}v(\lambda)t(\lambda')v(\lambda)^{-1}
 = u(\lambda')v(\lambda')v(\lambda')^{-1}
 = u(\lambda').\tag*{\qed}
 \]\renewcommand{\qed}{}
\end{proof}

\begin{Definition}
 For $\lambda'\in L'$, $\lambda\in L$, define
 \begin{align*}
& r'_k(\lambda'):= u(\lambda')(-\rho'_k)= v(\lambda')^{-1}(\lambda'_--\rho'_k),\\
& r_{k'}(\lambda):= u'(\lambda)(-\rho_{k'}) = v(\lambda)^{-1}(\lambda_--\rho_{k'}).
 \end{align*}
\end{Definition}

\begin{Lemma}[{\cite[Section~2.8.4]{macdonald}}] \label{sec:lem-pull-out-of-r}\leavevmode
\begin{enumerate}\itemsep=0pt
\item[$(i)$] If $u\in\Omega$ and $\lambda'\in L'$, then $r'_k(u\lambda')=ur'_k(\lambda')$.
\item[$(ii)$] If $i\in I_0$ with $s_i\lambda'\ne\lambda'$, then $r'_k(s_i\lambda')=s_ir'_k(\lambda')$.
\end{enumerate}
\end{Lemma}

\subsection{Partial order on the lattices}\label{sec:partial-order}
We shall now define partial orders on $L$, $L'$ that have finite
down-sets and thus allow for the application of the Gram--Schmidt algorithm.
These
partial orders will turn out to be intimately related to the Bruhat
ordering and the elements \smash{$(u(\lambda'))_{\lambda'\in L'}$} (resp.\
\smash{$(u'(\lambda))_{\lambda\in L}$}).

\begin{Definition}
 For $\lambda'\in L'$ and $\lambda\in L$ define $\overline{v}(\lambda')$ and
 $\overline{v}(\lambda)$ to be the shortest elements $w\in W_0$ mapping
 $\lambda'_+$ to $\lambda'$ and $\lambda_+$ to $\lambda$, respectively.
\end{Definition}

\begin{Definition}
 Let $\lambda',\mu'\in L'$, then define $\lambda'\le\mu'$ if
 \begin{enumerate}\itemsep=0pt
 \item[(i)] $\mu'-\lambda'\in\mathbb{N}_0 R^{\vee,+}\setminus\set{0}$ or
 \item[(ii)] $\mu'=\lambda'$ and $\overline{v}(\lambda')\le\overline{v}(\mu')$.
 \end{enumerate}
 Analogously on $L$.
\end{Definition}

\begin{Remark}
 We later refer to results from~\cite{stokmanInduced} although the authors
 of that paper use a~different definition for the order on~$L'$: they define
 $\lambda'\le\mu'$ iff $u(\lambda')\le u(\mu')$, i.e., they identify~$L'$
 with $W/W_0$ and take the Bruhat order on the shortest coset representatives.
 This definition is actually equivalent to our definition by
\cite[Proposition~5.21]{stokmanInduced} and \cite[Section~2.7.13]{macdonald}.
\end{Remark}

\begin{Lemma}\label{sec:lem-order-interval}
 Let $\lambda'\in L'_+$ and let $\mu'\in W_0\lambda'$. Let
 $\overline{v}(\mu')=s_{i_1}\cdots s_{i_p}$ be a reduced expression. Define
 $\mu'_r := s_{i_{r+1}}\cdots s_{i_p}\lambda'$, then
$
 \mu'=\mu'_0>\cdots>\mu'_p = \lambda'$.
\end{Lemma}
\begin{proof}
 Follows from \cite[Section~2.7.10\,(ii)]{macdonald} after some minor adjustments
 of notation.
\end{proof}

\subsection{Double-affine Hecke algebra}
For this subsection, we immediately give definitions for the (double-)affine
Hecke algebra in terms of generators and relations and thereby forego any mentions of the (double-)affine braid groups from \cite[Section~3]{macdonald}.
In order to still be able to use Macdonald's results, we prove that our definitions are equivalent.

We fix a field extension $K|\Q$ and a group homomorphism $q\colon \R\to K^\times$ such that $q(x)$ is transcendental over $\Q$ for all $x\ne0$
(e.g., $K=\Q(\R)$, the quotient field of the group algebra of~$\R$ over $\Q$). Write $A:= K[L]$ and $A':= K[L']$
for the group algebras of $L$, $L'$, whose monomials shall be written $e(\lambda)$, $e(\lambda')$. Furthermore, define $e(rc):=q(r)$, which provides
a definition for~$e(a)$ for any $a\in S$.

\begin{Definition}\label{def-tau}
 For a $W$-labelling $k\colon S\to\R$ and a root $a\in S$, define
 \[
 \tau_a := \tau_{a,k}:= q\qty(\frac{k(a)+k(2a)}{2}),\qquad
 \tilde{\tau}_a := \tilde{\tau}_{a,k} := q\qty(\frac{k(a)-k(2a)}{2}),
 \]
 where we take $k(2a)=0$ if $2a\not\in S$. In case $a=a_i$ for $i\in I$, we shall write
$\tau_i := \tau_{a_i}$, $ \tilde{\tau}_i := \tilde{\tau}_{a_i}$
 as a shorthand.
\end{Definition}

\begin{Remark}\label{sec:rem-taus}\leavevmode
\begin{enumerate}\itemsep=0pt
 \item[(i)] In the notation of Lemma~\ref{sec:lem-ks1},
 we can conclude
$
 \tau_i = \tau_{a'_i,k'}$.
 \item[(ii)] In the notation of Lemma~\ref{sec:lem-ks2}, we can conclude
$
 \tau_i = \tilde{\tau}_{a'_j,k'}$.
 \item[(iii)] In~\cite{macdonald}, the notation $\tau'_i$ is used instead of
 $\tilde{\tau}_i$, which can be easily confused as referring to the root system~$S'$ and the labelling~$k'$. To prevent such confusion, and in order to allow us
 to more easily talk about the $\tau$'s derived from~$k'$, we adopted the notation
 used here.
\end{enumerate}
\end{Remark}

\begin{Definition}
 Define the rational functions
 \begin{align*}
& \bm{b}(t,u;x):= \frac{t-t^{-1} + \bigl(u-u^{-1}\bigr)x}{1-x^2},\qquad
 \bm{c}(t,u;x) := \frac{tx - t^{-1}x^{-1} + u-u^{-1}}{x-x^{-1}}.
 \end{align*}
 Note that $\bm{c}\bigl(t^{-1},u^{-1};x^{-1}\bigr)=\bm{c}(t,u;x)$. Define
 furthermore the following functions:
 \[
 \bm{b}_a := \bm{b}_{a,k,S} := \bm{b}(\tau_{a,k},\tilde{\tau}_{a,k}, e(a))\in \operatorname{Frac}(A),\qquad
 \bm{b}'_{a'} := \bm{b}_{a',k',S'}\in\operatorname{Frac}(A')
 \]
 for $a\in S$, $a'\in S'$.
\end{Definition}

\begin{Definition}
 The \emph{affine Hecke algebra} $\mathfrak{H}$ is the associative $K$-algebra generated by
 $T(w)$ for $w\in W$ subject to
$
 T(v)T(w)=T(vw)$ if $ \ell(vw)=\ell(v)+\ell(w)$
 and
 \[
 (T(s_i)-\tau_i)\bigl(T(s_i)+\tau_i^{-1}\bigr) = 0
 \]
 for $i\in I$. For $i\in I$, write $T_i:= T(s_i)$. Then
 $\mathfrak{H}$ is in particular generated by $(T_i)_{i\in I}$ and~$(T(u))_{u\in\Omega}$. Since all elements of $\Omega$ have length 0,
 the latter family is a subgroup, and the $T_i$ satisfy the braid relations
 of the corresponding Coxeter system.
\end{Definition}

\begin{Lemma}\label{sec:lem-bernstein}
 There exists a group homomorphism $Y\colon L'\!\to\!\mathfrak{H}^\times$ with
 $Y^{\lambda'}\! = T(t(\lambda'))$ for~\smash{$\lambda'\in L'_+$}. Furthermore,
 $\mathfrak{H}$ is generated by \smash{$Y^{L'}$} and $(T_i)_{i\in I_0}$ subject to the
 following Bernstein--Lusztig--Zelevinsky relation:\footnote{Also referred to as Bernstein relation, Bernstein--Zelevinsky cross relation, or Lusztig's relation.}
 \begin{equation}\label{eq:y-bernstein}
 Y^{\lambda'} T_i - T_i Y^{s_i\lambda'} = \bm{b}'_{a_i'}\bigl(Y^{-1}\bigr)
 \bigl(Y^{\lambda'} - Y^{s_i\lambda'}\bigr)
 \end{equation}
 for $\lambda'\in L'$, $i\in I_0$,
 where \smash{$\bm{b}'_{a_i'}\bigl(Y^{-1}\bigr)$} means that we extend the group homomorphism
 $Y^{-1}\colon \lambda'\mapsto Y^{-\lambda'}$, extend it to $\operatorname{Frac}(A')$, and then apply it to \smash{$\bm{b}'_{a_i'}$}.
\end{Lemma}
\begin{proof}
 The existence of $Y$ follows from \cite[equation~(3.2.2)]{macdonald}, the
 relation, from \cite[Section~4.2.4]{macdonald}, and that $\mathfrak{H}$ is
 generated as claimed, from \cite[Section~4.2.7]{macdonald}.
\end{proof}

\begin{Lemma}[{\cite[Section~3.1.8]{macdonald}}]\label{sec:permute-simple-reflections}
 If $i,j\in I$ and $w\in W$ such that
 $ws_iw^{-1} = s_j$, then
$
 T(w)T_iT(w)^{-1} = T_j$.
\end{Lemma}

\begin{Definition}
 The \emph{double-affine Hecke algebra} $\tilde{\mathfrak{H}}$ is
 the associative $K$-algebra generated by~$\mathfrak{H}$ and elements $X^\lambda$ ($\lambda\in L$) such that
 \begin{enumerate}\itemsep=0pt
 \item[(i)] $\mathfrak{H}$ is a subalgebra and $X\colon L\to\tilde{\mathfrak{H}}^\times$
 is a group homomorphism. We extend $X$ to affine maps on $E$ whose gradient lies
 in $L$ by defining $X^c := q(1)$.
 \item[(ii)] The following Bernstein--Lusztig--Zelevinsky relation holds: Let $i\in I$ and $f\in F$ with $Df\in L$. Then
 \begin{equation}\label{eq:x-bernstein}
 T_i X^f - X^{s_if}T_i = \bm{b}_{a_i}(X)\bigl(X^f - X^{s_if}\bigr).
 \end{equation}
 \item[(iii)] $T(u) X^f T\bigl(u^{-1}\bigr) = X^{uf}$ for $u\in\Omega$.
 \end{enumerate}
\end{Definition}
\begin{Remark}
 In \cite[Section~4.7]{macdonald}, as well as other sources
 like \cite[Definition~2.1]{cherednikInduced}, $\tilde{\mathfrak{H}}$
 is defined with other relations, most of which are immediately a special case of
 the Bernstein presentation. Macdonald's use of potentially non-reduced affine root
 systems necessitates the imposition of another set of quadratic Hecke relations
 (that turn out to be only necessary for components of type $\bigl(C_n^\vee, C_n\bigr)$, i.e., for exactly these non-reduced root systems):
 we define
$\tilde{T}_i := X^{-a_i} T_i^{-1}$
 for~${i\in I}$. Then
 \begin{align*}
 \tilde{T}_i^{-1} - \tilde{\tau}_i^{-1} &=
 T_i X^{a_i} - \tilde{\tau}_i= X^{-a_i}T_i + \bm{b}_{a_i}(X)(X^{a_i}-X^{-a_i}) - \tilde{\tau}_i^{-1}\\
 &= X^{-a_i}T_i + \frac{\tau_i - \tau_i^{-1} + \bigl(\tilde{\tau}_i-\tilde{\tau}_i^{-1}\bigr)X^{a_i}}{1-X^{2a_i}} (X^{a_i}-X^{-a_i}) - \tilde{\tau}_i^{-1}\\
 &= X^{-a_i}T_i - X^{-a_i}\qty(\tau_i - \tau_i^{-1})
 - \tilde{\tau}_i\\
 &= X^{-a_i}\bigl(T_i - \tau_i + \tau_i^{-1}\bigr) - \tilde{\tau}_i= X^{-a_i} T_i^{-1} - \tilde{\tau}_i= \tilde{T}_i - \tilde{\tau}_i.
 \end{align*}

 This shows that the quadratic Hecke relations also hold for $\tilde{T}_i$, and therefore, that our definition for the double-affine
 Hecke algebra is indeed compatible with~\cite{cherednikInduced,macdonald}.
\end{Remark}

\begin{Lemma}[{\cite[Section~4.7.5]{macdonald}}]\label{sec:lem-basis-daha}
 The families
 \[
 \bigl(Y^{\lambda'}T(w)X^\mu\bigr)_{\lambda'\in L', w\in W_0, \mu\in L},\qquad
 \bigl(X^\mu T(w) Y^{\lambda'}\bigr)_{\mu\in L, w\in W_0, \lambda'\in L'}
 \]
 are $K$-bases of
 $\tilde{\mathfrak{H}}$ and hence give rise to Poincar\'e--Birkhoff--Witt-like
 decompositions
 \[
 \tilde{\mathfrak{H}}\cong A(X)\otimes\mathfrak{H}_0\otimes A'(Y)
 \cong A'(Y)\otimes\mathfrak{H}_0\otimes A(X),
 \]
 where $\mathfrak{H}_0\le\mathfrak{H}$ is the $($finite-dimensional$)$
 subalgebra generated by $(T_i)_{i\in I_0}$. These decompositions are
 as vector spaces and as left $($resp.\ right$)$ $A$-modules and as
 right $($resp.\ left$)$ $A'$-modules. In particular, we have
$
 \tilde{\mathfrak{H}}\cong A(X)\otimes\mathfrak{H} \cong \mathfrak{H}\otimes A(X)
$
 as left $($resp.\ right$)$ $A$-modules and as right $($resp.\ left$)$
 $\mathfrak{H}$-modules.
\end{Lemma}

\begin{Definition}
 Fix an involutive $\mathbb{Q}$-automorphism $*\colon K\to K$ mapping $q(x)$ to $q(-x)$.
 This exists because we required all $q(x)$ ($x\ne0$) to be transcendent.
\end{Definition}

\begin{Lemma}\label{sec:lem-existence-involution}
 $*$ extends to an involutive anti-automorphism of $\tilde{\mathfrak{H}}$ satisfying
 \[
 (X^\mu)^* = X^{-\mu},\qquad
 \bigl(Y^{\lambda'}\bigr)^* = Y^{-\lambda'},\qquad
 (T(w))^* = T(w)^{-1}
 \]
 for $\mu\in L$, $\lambda'\in L'$, $w\in W_0$.
\end{Lemma}
\begin{proof}
 Write $*$ also for the involutive automorphisms of $A$, $A'$ that are derived from
 the inverse map of $L$, $L'$. If $t^*=t^{-1}$, $u^*=u^{-1}$, $x^*=x^{-1}$, we have
 \begin{align*}
 \bm{b}(t,u;x)^* &= \frac{t^{-1} - t + \bigl(u^{-1}-u\bigr)x^{-1}}{1-x^{-2}}
 = \frac{x^2\bigl(t-t^{-1}\bigr) + x\bigl(u-u^{-1}\bigr)}{1-x^2}= \bm{b}(t,u;x) - t + t^{-1}.
 \end{align*}
 Consequently, it remains to show that the two Bernstein--Lusztig--Zelevinsky relations as well
 as the relations between $\Omega$ and $X^L$ are preserved by $*$.

 For $i\in I$ and $Df\in L$, we have
 \begin{align*}
 \bigl(T_i X^f - X^{s_if}T_i\bigr)^*
 &= X^{-f} T_i^{-1} - T_i^{-1} X^{-s_if}\\
 &= X^{-f} \bigl(T_i X^{s_if} - X^f T_i - \bigl(\tau_i-\tau_i^{-1}\bigr)\bigl(X^{s_if}- X^f\bigr)\bigr) X^{-s_if}\\
 &= X^{-f} \bigl(\bigl(\bm{b}_{a_i}(X) - \tau_i + \tau_i^{-1}\bigr)\bigl(X^{s_if} - X^f\bigr)\bigr)
 X^{-s_if}\\
 &= \bigl(\bm{b}_{a_i} - \tau_i + \tau_i^{-1}\bigr)(X)\bigl(X^{-f} - X^{-s_if}\bigr)\\
 &= \bm{b}_{a_i}(X)^* \bigl(X^f - X^{s_i f}\bigr)^*.
 \end{align*}
 A similar calculation shows the correct relations between
 $Y^{\lambda'}$ and $(T_i)_{i\in I_0}$, and for $\Omega$ note that $T(u)^*=T(u)^{-1}$ for $u\in\Omega$, so that
 \[
 \bigl(T(u) X^f T(u)^{-1}\bigr)^*
 = T(u) X^{-f} T(u)^{-1}
 = X^{u(-f)}
 = \bigl(X^{u(f)}\bigr)^*
 \]
 since $\Omega$ acts linearly on $F$.
\end{proof}

\subsection{Basic representation and inner product}
We shall now describe the basic representation of $\tilde{\mathfrak{H}}$.
Contrary to the approach of~\cite{macdonald} but in anticipation of Section~\ref{sec:jasper}, we shall describe the basic representation as an induced module.

\begin{Lemma}\label{sec:lem-trivial-character}
 There is a unique ring homomorphism $\phi\colon \mathfrak{H}\to K$ satisfying
 $\phi(T_i) = \tau_{a_i',k'}=\tau_i$ for $i\in I_0$ and
$
 \phi\bigl(Y^{\lambda'}\bigr) = q(\langle \lambda', \rho_{k'}\rangle)$.
\end{Lemma}
\begin{proof}
 By Lemma~\ref{sec:lem-bernstein}, $\mathfrak{H}$ is generated by $(T_i)_{i\in I_0}$ and $Y^{L'}$ satisfying~\eqref{eq:y-bernstein}.
 Thus, any ring homomorphism is uniquely defined by what it maps $Y^{L'}$ and the $T_i$ to. This shows uniqueness.

 For existence, we need to show that the given data for $\phi$
 satisfies the relations of $\mathfrak{H}$. The relations
 among the $T_i$ are the Hecke relations (which are trivially satisfied) and the braid relations of the corresponding Coxeter
 system. Let $(m_{ij})_{i,j\in I}$ be the Coxeter matrix and let~${i,j\in I_0}$ with~${i\ne j}$. If $m_{ij}$ is even, then the
 braid relation is satisfied because $\tau_i$ and $\tau_j$ commute with each other. If $m_{ij}$ is odd, the
 roots $a_i'$ and $a_j'$ lie in the same $W_{S'}$-orbit and therefore,
 $\tau_i=\tau_j$. This shows that the braid relations
 are also satisfied for this case.

 It thus remains to show the Bernstein--Lusztig--Zelevinsky relation, so let $\lambda'\in L'$
 and $i\in I_0$. From Lemma~\ref{sec:lem-rho-inner-prod}, we know
 that $\langle\alpha_i^\vee,\rho_{k'}\rangle=k'(a_i')$. Furthermore,
 $\tilde{\tau}_{a'_i,k'}q(-k'(a_i'))=\tau_{a'_i,k'}^{-1}$, so that
 \begin{align*}
 \frac{\tau_{a_i',k'} - \tau_{a_i',k'}^{-1} + \bigl(\tilde{\tau}_{a'_i,k'} - \tilde{\tau}_{a'_i,k'}^{-1}\bigr)q(-k'(a_i'))}{1-q(-2k'(a'_i))}
& = \frac{\tau_{a_i',k'} - \tau_{a'_i,k'}^{-1}
 + \tau_{a'_i,k'}^{-1} - \tau_{a'_i,k'}q(-2k'(a_i'))}{1-q(-2k'(a'_i))}\\
& = \tau_{a'_i, k'} = \tau_i
 \end{align*}
 by Remark~\ref{sec:rem-taus}.
 As a consequence, we have
 \begin{gather*}
 \frac{\tau_{a_i',k'} - \tau_{a_i',k'}^{-1} + \bigl(\tilde{\tau}_{a'_i,k'} - \tilde{\tau}_{a'_i,k'}^{-1}\bigr)q(-k'(a_i'))}{1-q(-2k'(a'_i))} (q(\langle\lambda',\rho_{k'}\rangle)-q(\langle s_i\lambda',\rho_{k'}\rangle))\\
\qquad = q(\langle\lambda',\rho_{k'}\rangle)\tau_i
 - \tau_iq(\langle s_i\lambda',\rho_{k'}\rangle),
 \end{gather*}
 which shows that our data satisfies~\eqref{eq:y-bernstein}, and hence gives rise to a ring homomorphism \linebreak $\mathfrak{H}\to K$.
\end{proof}

Via the ring homomorphism $\phi$ from Lemma~\ref{sec:lem-trivial-character}, we can view $K$ as a $\mathfrak{H}$-module and then consider
the $\tilde{\mathfrak{H}}$-module $M:=\tilde{\mathfrak{H}}\otimes_{\mathfrak{H}}\otimes K$.

\begin{Lemma}\label{sec:lem-phi-values}
 For $i\in I$, we have $\phi(T_i)=\tau_i$ $($in particular, also for
 $i\not\in I_0)$, and for $u\in\Omega$, we have
$\phi(T(u))=1$.
\end{Lemma}
\begin{proof}
 For $i\in I_0$, we have $\phi(T_i)=\tau_i$ by Remark~\ref{sec:rem-taus}\,(i). Otherwise,
 by \cite[equation~(3.3.5)]{macdonald}, if $w\in W_0$ and
 $w\alpha_j\in R$ is the highest
 root in the same irreducible component as $a_i$, we have
 \[
 T_i = T(w) Y^{\alpha_j^\vee} T_j^{-1} T(w)^{-1}.
 \]
 Applying $\phi$, we obtain
 \[
 \phi(T_i) = \tau_j^{-1} q\bigl(\big\langle\alpha_j^\vee,\rho_{k'}\big\rangle\bigr)
 = \tau_{a'_j,k'}^{-1} q(k'(a'_j))
 \]
 by Lemma~\ref{sec:lem-rho-inner-prod}. Observe that by the proof of
 Lemma~\ref{sec:lem-trivial-character}, this equals
 \smash{$\tilde{\tau}_{a'_j,k'}$}, which equals $\tau_i$ by Remark~\ref{sec:rem-taus}\,(i).

 For $u$, note that $\phi$ restricted to $T(\Omega)$ is a character of
 $\Omega$. Since we can write $T(u)$ as a product of $T_i$'s and elements of $Y^{L'}$, and since $\phi$ is defined in terms of $q$,
 we have $\phi(T(u))=q(x)$ for some $x\in\R$. Since $\Omega$ is a finite group,
 $u$ has finite order, say $n$. Then $\phi(T(u)^n)=\phi(T(u^n))=1=q(nx)$.
 As $q$ is injective, we have $nx=0$ and hence $x=0$.
\end{proof}

\begin{Remark}
 Evidently, it would have been easier to define $\phi$ using the
 relations from Lem\-ma~\ref{sec:lem-phi-values}. The approach chosen here
 makes a clearer connection to the content of Sections~\ref{sec:parabolically-induced} and~\ref{sec:spherical-vec}.
\end{Remark}

\begin{Lemma}\label{sec:lem-induced-rep-iso}
 The map $e(\mu)\mapsto X^\mu\otimes 1$ defines an isomorphism
 $A\cong M$ of $K$-vector spaces and of left $A$-modules.
\end{Lemma}
\begin{proof}
 By Lemma~\ref{sec:lem-basis-daha}, the elements $(X^\mu)_{\mu\in L}$ are
 a $\mathfrak{H}$-basis (with respect to the action from the right),
 and thus $(X^\mu\otimes 1)_{\mu\in L}\subset M$ are $K$-linearly independent
 and a $K$-generating system. Since our described map maps a $K$-basis
 of $A$ to a $K$-basis of $M$, it is a $K$-linear isomorphism~${A\cong M}$. Since it also respects the $A$-module structures on $A$ and
 $M$, it is even an $A$-linear isomorphism.\looseness=1
\end{proof}

\begin{Definition}
 $A$, viewed as an $\tilde{\mathfrak{H}}$-module using the isomorphism from
 Lemma~\ref{sec:lem-induced-rep-iso} is called the \emph{basic representation}.
\end{Definition}

\begin{Lemma}\label{sec:lem-concrete-basic-rep}
 Let $f\in A$, then we have
$
 T_i f = \tau_i s_i f + \bm{b}_{a_i} (f - s_i f)$, $
 T(u) f = uf$, $
 X^\mu f = e(\mu)f
$
 for $i\in I$, $u\in\Omega$, $\mu\in L$.
\end{Lemma}
\begin{proof}
 By~\eqref{eq:x-bernstein} and linear extension, we have
 \[
 T_i f(X)\otimes 1 = (s_if)(X)\otimes \phi(T_i) + \bm{b}_{a_i}(X)(f-s_if)(X)\otimes 1.
 \]
 By Lemma~\ref{sec:lem-phi-values}, we have $\phi(T_i)=\tau_i$, whence
 we have the first claim. For the second claim, let $u\in \Omega$, then
 \[
 T(u)f(X)\otimes 1 = (uf)(X)\otimes\phi(u)
 = (uf)(X)\otimes 1
 \]
 by Lemma~\ref{sec:lem-phi-values}. Lastly, the third claim follows from the fact that our identification
 is $A$-linear.
\end{proof}

We can equip the basic representation with an inner product and make it
unitary.

\begin{Definition}\label{sec:def-deltas-nabla}
 Let $\Delta=\Delta_{k,S}$ be
 an appropriately defined infinite product
 \[
 \prod_{a\in S_1^+} \tau_a^{-1}\bm{c}(\tau_a,\tilde{\tau}_a, e(a))^{-1}
 \]
 as in \cite[Section~5.1.7\,ff]{macdonald}.
 Define furthermore
 \begin{align*}
 \Delta_0:= \Delta_{0,S}:= \prod_{\substack{a\in S_1^+\\Da=0}}
 \tau_a^{-1}\bm{c}(\tau_a,\tilde{\tau}_a,e(-a))^{-1},\qquad
 \nabla:= \nabla_S :=
 \Delta\Delta_0.
 \end{align*}
\end{Definition}

\begin{Lemma}\label{lem-weight-distributions}\quad
 \begin{enumerate}\itemsep=0pt
 \item[$(i)$] There are distributions $\Delta_1,\nabla_1\in K[[L]]$ that are proportional to
 $\Delta$, $\nabla$ and have $1$ as their constant term;
 \item[$(ii)$] $\nabla$ is $W_0$-symmetric;
 \item[$(iii)$] we have
$
 \sum_{w\in W_0} w\Delta_0^{-1} = W_0\bigl(\tau^2\bigr)
$
$($already using notation from Example~{\rm \ref{sec:ex-labellings}\,(i), (ii))}.
 \end{enumerate}
\end{Lemma}
\begin{proof}
(i) Follows by \cite[Section~5.1.10]{macdonald}
 and the fact that $\Delta_0\in K(L)$.
(ii) Follows by \cite[Section~5.1.27]{macdonald}.
(iii) Follows by \cite[Section~5.1.36]{macdonald}.
\end{proof}

Recall that we can use any distribution $h\in K[[L]]$
to define a Hermitean bilinear form on~$A$ as follows
$
 (f,g) \mapsto \ct(fg^* h)
$
and a symmetric bilinear form on~$A$ as follows
$
 (f,g)\mapsto \ct(f\overline{g}h)$.
Here $\ct\colon K[[L]]\to K$ maps a distribution to its constant term, and $*$ is the involution of $A$ extending $*$ on $K$
(i.e., mapping $q(a)\mapsto q(-a)$ for $a\in\Q$) and mapping
$e(\mu)\mapsto e(-\mu)$ for $\mu\in L$.
Furthermore, $\overline{\cdot}$ is the $K$-linear
involution mapping $e(\mu)\mapsto e(-\mu)$ for $\mu\in L$ (but preserving $q(a)$ for $a\in\Q$).

We thus define our inner products as follows:
\begin{Definition}\label{def-inner-product}
 For $f,g\in A$, define
$
 (f,g):= (f,g)_{k,S} := \ct(fg^*\Delta_{k,S})$.
 Furthermore, define
\smash{$
 (f,g)_1 := \frac{(f,g)}{(1,1)}
 = \ct(fg^*\Delta_{k,S,1})\in K$}.
\end{Definition}

By \cite[Section~5.1.20]{macdonald}, the equation $(f,f)=0$ implies $f=0$.
Consequently, we can apply the Gram--Schmidt procedure
to $(e(\mu))_{\mu\in L}$, with the monomial order from Section~\ref{sec:partial-order}.
\begin{Theorem}[{\cite[Section~5.2.1ff]{macdonald}}]\label{sec:thm-E}
 There is a unique family $(E_\lambda)_{\lambda\in L}\subset A$ satisfying
 \begin{enumerate}\itemsep=0pt
 \item[$(i)$] $E_{\lambda} = e(\lambda) + \text{\rm l.o.t.}$ {\rm(}here, ``lower'' refers
 to the order $\le$ on $L${\rm)};
 \item[$(ii)$] $f(Y) E_\lambda = f(-r_{k'}(\lambda)) E_\lambda$;
 \item[$(iii)$] $(E_\lambda, E_\mu)=0$ for $\lambda\ne\mu$,
 \end{enumerate}
 called the non-symmetric Macdonald polynomials.
 In particular, the span of $(E_\mu)_{\mu\in W_0\lambda}$ is
 an $\mathfrak{H}_0$-module.
\end{Theorem}

That last remark of the $(E_\mu)_{\mu\in W_0\lambda}$ being an
$\mathfrak{H}_0$-module can be made more precise:
\begin{Lemma}[{\cite[Section~5.4.3f]{macdonald}}]\label{sec:lem-Ti-action-E}
 Let $\lambda\in I$ and $i\in I_0$.
 \begin{enumerate}\itemsep=0pt
 \item[$(i)$] If $\langle\lambda, a'_i\rangle=0$, so that $s_i\lambda=\lambda$,
 we have $T_iE_\lambda = \tau_i E_\lambda$.
 \item[$(ii)$] If $\langle\lambda, a'_i\rangle>0$, so that $s_i\lambda>\lambda$,
 we have $T_iE_\lambda = \tau_i^{-1} E_{s_i\lambda} + \bm{b}'_{a'_i}(r_{k'}(\lambda))E_\lambda$.
 \end{enumerate}
\end{Lemma}

\section{Adapting Macdonald's formalism to parabolic subgroups}\label{sec:macdo-parabolic}
Fix a subset $J\subset I_0$ and write $W_J\le W_0$ and $\mathfrak{H}_J\le\mathfrak{H}_0$ for
the subgroup/-algebra generated by $s_j$ and $T_j$ ($j\in J$), respectively. Let
$A_J := A^{W_J}$, $A'_J:= (A')^{W_J}$, and let $W_0^J$ be the set of shortest representatives of elements of $W_0/W_J$. Let
$w_J$, $w_0$ be the longest elements of~$W_J$,~$W_0$, respectively. We shall now
see how the formalism from~\cite{macdonald} can be used in a $W_J$-invariant
context.

\subsection{Poincar\'e polynomials}
We begin by recalling some facts about (much generalised) Poincar\'e
polynomials and series, as it turns out, much of our theory requires them.
Throughout this subsection, let $(\mathcal{W},\mathcal{S})$ be a finite Coxeter group and
$\mathcal{A}$ a ring.

\begin{Definition}
 A function
 $\tau\colon \mathcal{W}\to \mathcal{A}^\times$ is called
 an ($\mathcal{A}$-valued) \emph{multiplicative labelling} of $\mathcal{W}$ if $\ell(vw)=\ell(v)+\ell(w)$ implies that
 $\tau_{vw}=\tau_v\tau_w$.

 For any multiplicative labelling $\tau$ and a subset $X\subset \mathcal{W}$, define
$
 X(\tau) := \sum_{w\in X} \tau_w$.
 This is the \emph{Poincar\'e polynomial} of the subset $X$ with respect to $\tau$.
\end{Definition}

\begin{Remark}
 The usual meaning of ``Poincar\'e polynomial'' can be recovered by
 considering~${A=\mathbb{Z}[t]}$ and \smash{$\tau_w := t^{\ell(w)}$}.
\end{Remark}

\begin{Lemma}\label{sec:lem-construct-labelling}
 Let $(\tau_s)_{s\in \mathcal{S}}\subset \mathcal{A}^\times$ satisfy the braid relations of $(\mathcal{W},\mathcal{S})$.
Then the map $\tau\colon \mathcal{S}\to \mathcal{A}^\times$, $s\mapsto \tau_s$ can be extended
 uniquely to a multiplicative labelling of~$\mathcal{W}$.
\end{Lemma}
\begin{proof}
 In analogy with \cite[Section~3.1]{macdonald}, the Artin braid group
 $\mathfrak{B}$ of $\mathcal{W}$ is generated by $(T(w))_{w\in \mathcal{W}}$ subject to
$
 \ell(vw)=\ell(v)+\ell(w)\Rightarrow T(v)T(w)=T(vw)$.
 It follows therefore that the multiplicative labellings (with values in
 $\mathcal{A}$) of $\mathcal{W}$ are exactly the representations of
 $\mathfrak{B}$ (in~$\mathcal{A}^\times$). By
 \cite[Section~3.1.6]{macdonald}, the braid group can also be
 presented in terms of the generators~$(T_s)_{s\in \mathcal{S}}$ and the braid
 relations. Thus a representation can be specified by providing
 elements~$(\tau_s)_{s\in\mathcal{S}}$ satisfying the braid relations.
\end{proof}

\begin{Example}\label{sec:ex-labellings}\qquad
\begin{enumerate}\itemsep=0pt
 \item[(i)] If $\tau$ is a multiplicative labelling that maps to a commutative algebra,then
 any integer power of $\tau$ is a multiplicative labelling as well.
 \item[(ii)]
 Let $k$ be a $W$-labelling of $S$. For $j\in J$, take $\tau_j := \tau_{a_j,k}$ as in Definition~\ref{def-tau}. These elements satisfy the braid relations, so that
 $\tau=\tau_k$ is a multiplicative labelling of $W_J$ (after
 appropriate restriction).
 \item[(iii)] Let in addition $\epsilon\colon W_J\to\C^\times$ be a multiplicative character. Then
 define
 \[
 \tau^{(\epsilon)}_j := \tau^{(\epsilon)}_{j,k} :=
 \begin{cases}
 \tau_{j,k}, & \epsilon(s_j)=1,\\
 -\tau_{j,k}^{-1}, & \epsilon(s_j)=-1,
 \end{cases}
 \]
 the \emph{$q$-deformation of $\epsilon$}. This also gives rise to
 a multiplicative labelling \smash{$\tau^{(\epsilon)}_k$} of $W_J$, and even
 to a ring
 homomorphism $\mathfrak{H}_J\to K$.
 The name $q$-deformation comes from the fact that $\mathfrak{H}_J$
 is a deformation of $K[W_J]$, and that $\epsilon$ can be extended to
 a ring homomorphism $K[W_J]\to K$ that can be obtained from~\smash{$\tau^{(\epsilon)}_k$} as a limit case.
 \item[(iv)] Take $A=\mathfrak{H}_J$, then the $(T_j)_{j\in J}$ satisfy the braid
 relations, giving rise to a multiplicative labelling~$T$.
 \item[(v)] Similarly, the \smash{$(\tau_{j,k} T_j)_{j\in J}$} also satisfy the braid
 relations, as do \smash{$\bigl(\tau^{(\epsilon)}_{j,k} T_j\bigr)_{j\in J}$}, yielding multiplicative labellings $\tau_k T$ and \smash{$\tau^{(\epsilon)}_k T$}.
\end{enumerate}
\end{Example}

\begin{Lemma}\label{sec:lem-decompose-poincare}
 Let $\mathcal{J}\subset\mathcal{S}$ a subset. Let $\tau$ be a
 multiplicative labelling on $\mathcal{W}$, then we have the following identity of Poincar\'e polynomials:
$
 \mathcal{W}(\tau) = \mathcal{W}^{\mathcal{J}}(\tau) \mathcal{W}_{\mathcal{J}}(\tau)$,
 where $\mathcal{W}_{\mathcal{J}}$ is the subgroup generated by $\mathcal{J}$, and $\mathcal{W}^{\mathcal{J}}$ is the set of shortest coset representatives of $\mathcal{W}/\mathcal{W}_{\mathcal{J}}$.
\end{Lemma}
\begin{proof}
 For the case \smash{$\tau_w=t^{\ell(w)}$} or more generally, $\mathcal{A}$
 commutative, this result is well known, e.g., \cite[Section~5.12]{humphreysRefl}, but it even holds for $\mathcal{A}$ non-commutative.

 For every $w\in\mathcal{W}$, there is a unique decomposition $vw'$ with
 $v\in \mathcal{W}^{\mathcal{J}}$ and $w'\in \mathcal{W}_{\mathcal{J}}$ such that~${\ell(w)=\ell(v)+\ell(w')}$. Therefore,
 \begin{align*}
& \mathcal{W}(\tau)= \sum_{v\in \mathcal{W}^{\mathcal{J}},\, w\in \mathcal{W}_{\mathcal{J}}}
 \tau_{vw} = \sum_{v\in \mathcal{W}^{\mathcal{J}},\, w\in \mathcal{W}_{\mathcal{J}}}
 \tau_v \tau_w= \sum_{v\in \mathcal{W}^{\mathcal{J}}} \tau_v \sum_{w\in \mathcal{W}_{\mathcal{J}}} \tau_w
 = \mathcal{W}^{\mathcal{J}}(\tau) \mathcal{W}_{\mathcal{J}}(\tau).\tag*{\qed}
 \end{align*} \renewcommand{\qed}{}
\end{proof}

\begin{Lemma}\label{sec:lem-inverse-labelling}
 For any multiplicative labelling $\tau\colon \mathcal{W}\to \mathcal{A}^\times$, the map $\tau^{-1}\colon w\mapsto \tau_w^{-1}$ is an
``anti-multiplicative'' labelling: we have $\tau_{vw}^{-1} = \tau_w^{-1}\tau_v^{-1}$ whenever $\ell(vw)=\ell(v)+\ell(w)$.
 We then have
$
 \mathcal{W}(\tau^{-1}) = \tau_{w_0}^{-1} \mathcal{W}(\tau)
$
$($where $\mathcal{W}(\tau^{-1})$ is defined as if $\tau^{-1}$ were a multiplicative labelling,
 and where $w_0$ is the longest element of $\mathcal{W})$.
\end{Lemma}
\begin{proof}
 Let $w\in \mathcal{W}$. By \cite[Corollary~1.8]{humphreysRefl}, we have
 $\ell(w_0w)=\ell(w_0)-\ell(w)$, i.e., $\ell(w_0w)+\ell(w)=\ell(w_0)$.
 Consequently, we have $\tau_{w_0w}\tau_w=\tau_{w_0}$, or
 equivalently
 $\tau_w^{-1}=\tau_{w_0}^{-1}\tau_{w_0w}$.
 Thus,
 \[
 \mathcal{W}(\tau^{-1}) = \sum_{w\in\mathcal{W}} \tau_w^{-1}
 = \tau_{w_0}^{-1}\sum_{w\in\mathcal{W}} \tau_{w_0w}
 = \tau_{w_0}^{-1}\sum_{w\in\mathcal{W}} \tau_w
 = \tau_{w_0}^{-1} \mathcal{W}(\tau).\tag*{\qed}
 \]\renewcommand{\qed}{}
\end{proof}

\begin{Remark}
 All results of this subsection except for Lemma~\ref{sec:lem-inverse-labelling} also work in case $(\mathcal{W},\mathcal{S})$ is an infinite
 Coxeter group, provided sufficient convergence of the series, e.g., $\mathcal{A}$ being filtered and $\tau$ mapping to $\mathcal{A}^+$.
 In this case, we speak not of Poincar\'e polynomials, but of \emph{Poincar\'e series}.
\end{Remark}

\subsection{Symmetrisers}
From this general finite Coxeter group interlude, we now return to our
main programme of affine root systems $S$, $S'$, their extended affine Weyl
groups $W$, $W'$, and the parabolic subgroup ${W_J\le W_0}$. Using the
language of (generalised) Poincar\'e polynomials, we now define symmetrisers
in a way generally inspired by \cite[Section~5.5]{macdonald}.

\begin{Definition}
 Let $\epsilon\colon W_J\to\C^\times$ be a multiplicative character of $W_J$.
 Define
 \[
 U^{(\epsilon)}_J := U^{(\epsilon)}_{J,k} := \bigl(\tau^{(\epsilon)}_{w_J,k}\bigr)^{-1} W_J\bigl(\tau^{(\epsilon)}_k T\bigr),
 \]
 the \emph{$\epsilon$-symmetriser} for $W_J$. Here, we use Example~\ref{sec:ex-labellings}\,(v).
\end{Definition}
We furthermore adopt the convention that for all notation that has an exponent $(\epsilon)$, we interpret the absence of that exponent to mean
the trivial character.

\begin{Lemma}\label{sec:lem-symmetriser-symmetric}
 Let $j\in J$, then
 \[
 \bigl(T_j-\tau^{(\epsilon)}_j\bigr)U^{(\epsilon)}_J =
 U^{(\epsilon)}_J \bigl(T_j-\tau^{(\epsilon)}_j\bigr) = 0.
 \]
\end{Lemma}
\begin{proof}
 Follows from \cite[Section~5.5.9]{macdonald} applied to $\mathfrak{H}_J$.
\end{proof}

This shows that for any representation of $\mathfrak{H}_J$, the application of \smash{$U^{(\epsilon)}_J$}
produces elements that transform according to the character $\tau^{(\epsilon)}$. In fact, \smash{$U^{(\epsilon)}_J$}
is (a scalar multiple of) the projection onto the corresponding isotypic component:
\begin{Corollary}\label{sec:cor-makes-symmetric}
 Let $V$ be a representation of $\mathfrak{H}_J$, then
 \[
 \big\{v\in V\mid \forall j\in J\colon T_j v = \tau^{(\epsilon)}_j v\big\}
 = U^{(\epsilon)}_J V.
 \]
\end{Corollary}
\begin{proof}
``$\supset$'': Follows from Lemma~\ref{sec:lem-symmetriser-symmetric}.

``$\subset$'': Let $v\in V$ have the desired transformation behaviour, then we have \smash{$\tau^{(\epsilon)}_w T(w) v = \tau^{(\epsilon)2}_w v$} for all $w$, hence
 \[
 U^{(\epsilon)}_J v
 = \tau_{w_J}^{-1} W_J\bigl(\tau^{(\epsilon)}T\bigr) v
 = \frac{W_J(\tau^{(\epsilon)2})}{\tau_{w_J}}v.
 \]
 Multiplying by an appropriate scalar shows that $v$ lies
 in $U^{(\epsilon)}_JV.$
\end{proof}

An immediate consequence of Lemma~\ref{sec:lem-decompose-poincare} is
the following statement about symmetrisers.

\begin{Corollary}\label{sec:lem-decomposition-symmetriser}
 Let $J'\subset J$ be a subset, then
 \[
 U^{(\epsilon)}_J = \bigl(\tau^{(\epsilon)}_{w_Jw_{J'}}\bigr)^{-1}\sum_{v\in W_J^{J'}} \tau^{(\epsilon)}_v
 T(v) U^{(\epsilon)}_{J'}
 = \bigl(\tau^{(\epsilon)}_{w_Jw_{J'}}\bigr)^{-1} W_J^{J'}\bigl(\tau^{(\epsilon)}T\bigr)
 U^{(\epsilon)}_{J'}.
 \]
\end{Corollary}

We conclude with some further properties of symmetrisers.
\begin{Proposition}\label{sec:prop-symmetriser-properties}
\qquad
 \begin{enumerate}\itemsep=0pt
 \item[$(i)$] \smash{$\bigl(U^{(\epsilon)}_J\bigr)^2=\frac{W_J(\tau^{(\epsilon)2}_k)}{\tau^{(\epsilon)}_{w_J}} U^{(\epsilon)}_J$}.\vspace{2mm}

 \item[$(ii)$] \smash{$\bigl(U^{(\epsilon)}_J\bigr)^*=U^{(\epsilon)}_J$}.
 \item[$(iii)$] For $f,g\in A$, we have
 \[
 \bigl(U^{(\epsilon)}_J f, U^{(\epsilon)}_Jg\bigr)
 = \frac{W_J\bigl(\tau^{(\epsilon)2}_k\bigr)}{\tau^{(\epsilon)}_{w_J}}
 \bigl(U^{(\epsilon)}_J f, g\bigr).
 \]
 \end{enumerate}
\end{Proposition}
\begin{proof} \cite[Section~5.5.17\,(ii)]{macdonald} applied
 to $\mathfrak{H}_J$.
 \cite[Section~5.5.17\,(iii)]{macdonald} applied to
 $\mathfrak{H}_J$.
Using the first two statements, we get
 \[
 \bigl(U^{(\epsilon)}_J f, U^{(\epsilon)}_J g\bigr)
 = \bigl(\bigl(U^{(\epsilon)}_J\bigr)^2 f, g\bigr)
 = \frac{W_J\bigl(\tau^{(\epsilon)2}_k\bigr)}{\tau^{(\epsilon)}_{w_J}}
 \bigl(U^{(\epsilon)}_J f, g\bigr).\tag*{\qed}
 \]
\renewcommand{\qed}{}
\end{proof}

\subsection[W\_J-orbits]{$\boldsymbol{ W_J}$-orbits}
In this subsection, we set out for a suitable way of labelling $W_J$-orbits within
$L$, and examine what~$\le$ looks like when restricted to a $W_J$-orbit.
\begin{Definition}
 Let
 \[
 L_{+,J} := \set{\lambda\in L\where\forall j\in J\colon \langle \lambda, a_j'\rangle\ge0} \subset L
 \]
 be the set of \emph{$J$-dominant} elements.
\end{Definition}

A useful property of $L_+$ is that all of its stabilisers are parabolic.
It turns out that $L_{+,J}$ possesses the same property.
\begin{Lemma}\label{sec:lem-J-stabiliser-parabolic}
 Let $\lambda\in L_{+,J}$, then the group
 \[
 W_{J,\lambda} := \set{w\in W_J\where w\lambda=\lambda}
 = W_\lambda\cap W_J
 \]
 is a parabolic subgroup of $W_J$.
\end{Lemma}
\begin{proof}
 Let $W'$ be the subgroup of $W_{J,\lambda}$ generated by its simple reflections. We now show
 that~${W_{J,\lambda}\subset W'}$. Let $w\in W_{J,\lambda}$; we proceed by induction in $p=\ell(w)$:

``$p=0$'': Then $w=1$, so $w\in W'$.

``$p-1\to p$'': Let $\ell(w)=p$. Since $1\ne w$, there is $j\in J$ such that
 $wa'_j<0$. Then we have $\ell(ws_j)<\ell(w)$. Furthermore, we have
 \[
 0 \le \langle \lambda, a'_j\rangle
 = \langle w\lambda, wa'_j\rangle
 = \langle \lambda, wa'_j\rangle \le 0
 \]
 since $\lambda\in L_{+,J}$ and since $wa'_j<0$ is a negative linear combination of
 $a_i$ ($i\in J$). Thus, $s_j\lambda=\lambda$.
 As a consequence, $ws_j\lambda = w\lambda=\lambda$. Since $\ell(ws_j)=p-1$,
 we can apply the induction hypothesis to find that $ws_j\in W'$, as is $s_j$, and
 hence also $w\in W'$.
\end{proof}

An immediate corollary to Lemma~\ref{sec:lem-J-stabiliser-parabolic} concerns the decomposition of
$W_J$-symmetrisers from Corollary~\ref{sec:lem-decomposition-symmetriser}.
\begin{Corollary}\label{sec:cor-decomposition-symmetriser}
 Let $\lambda\in L_{+,J}$ such that $\epsilon(W_{J,\lambda})=\set{1}$, then
 \[
 U^{(\epsilon)}_J = \bigl(\tau^{(\epsilon)}_{w_Jw_{J,\lambda}}\bigr)^{-1}\sum_{v\in W_J^\lambda} \tau^{(\epsilon)}_v
 T(v) U_{J'}
 \]
 where $J'=\{j\in J\mid \langle\lambda,a'_j\rangle=0\}$; in other words, where $W_{J'}=W_{J,\lambda}$.
\end{Corollary}

In order to prove that $L_{+,J}$ is a fundamental domain for $W_J$, it is useful
to know how to produce (the) elements of $L_{+,J}$.
\begin{Lemma}\label{sec:lem-obtain-J-dominant}
 Let $v\in W_0^J$ and $\lambda\in L_+$, then $v^{-1}\lambda\in L_{+,J}$.
\end{Lemma}
\begin{proof}
 Let $j\in J$. Since $v$ is the shortest element of its right coset, we have $\ell(vs_j)>\ell(v)$ and hence
 $va'_j>0$ by Lemma~\ref{sec:lem-length-change}. Since $\lambda$ is dominant, we therefore have
 \[
 \big\langle v^{-1}\lambda, a'_j\big\rangle
 = \langle \lambda, va'_j\rangle\ge0.\tag*{\qed}
 \]\renewcommand{\qed}{}
\end{proof}

\begin{Proposition}\label{sec:prop-J-fundamental-domain}
 The set $L_{+,J}$ is a fundamental domain for the action of $W_J$ on $L$.
\end{Proposition}
\begin{proof}
 Let $\lambda\in L$, we show that there is a unique element in $L_{+,J}\cap W_J\lambda$.

{\it Uniqueness.} Let $\lambda\in L_{+,J}$, we show that
$
 \set{w\in W_J\where w\lambda\in L_{+,J}}\subset W_{J,\lambda}$,
 which we do by induction on length $p$.

``$p=0$'': Evidently true.

``$p-1\to p$'': Let $\ell(w)=p$ and $w\lambda\in L_{+,J}$. Then there is $j\in J$ such
 that $\ell(s_jw)=p-1$. By Lemma~\ref{sec:lem-length-change}, this
 implies that $w^{-1}a'_j<0$ is a negative
 linear combination of $a'_i$ ($i\in J$), and hence that
 $\big\langle \lambda, w^{-1}a'_j\big\rangle \le 0$. Thus we have
$
 0 \ge \big\langle \lambda, w^{-1}a'_j\big\rangle = \langle w\lambda, a'_j\rangle \ge 0$,
 which implies that $s_j w\lambda = w\lambda\in L_{+,J}$. Since $\ell(s_jw)=p-1$, we can apply the
 induction hypothesis to find that $s_jw\in W_{J,\lambda}$ and in particular $s_jw\lambda=\lambda$, which
 then also shows that $w\lambda=\lambda$. As a consequence, every $W_J$-orbit contains at most one element
 in $L_{+,J}$.

{\it Existence.} We know that $L_+$ is a fundamental domain for $W_0$, so let
 $w\in W_0$ such that $\mu=w\lambda\in L_+$. Decompose $w=vw'$ for $v\in W_0^J$, $w'\in W_J$. Then
$
 w'\lambda = v^{-1}\mu\in L_{+,J}
$
 by Lemma~\ref{sec:lem-obtain-J-dominant}.
\end{proof}

\begin{Remark}
 The reader might be tempted to prove Proposition~\ref{sec:prop-J-fundamental-domain} using the fact that $W_J$ is the Weyl group of
 a~finite root system whose dominant weights are always a~fundamental
 domain for $W_J$'s action. However, this does not work here since
 we're considering $W_J$'s action on a~much larger lattice $L$.
\end{Remark}

We now want to understand the order restricted to a $W_J$-orbit.
This requires a more specialised version of $\overline{v}$.

\begin{Definition}
 Let $\lambda\in L_{+,J}$ and let $\mu\in W_J\lambda$. Write
 $\overline{v}_J(\mu)$ for the shortest element $w\in W_J$ satisfying $w\lambda=\mu$.
\end{Definition}

It turns out that this notion is related to $\overline{v}$ from Section~\ref{sec:partial-order}
via the shortest coset representatives decomposition.

\begin{Lemma}\label{sec:lem-decompose-vs}
 Let $\lambda\in L_{+,J}$ and $\mu\in W_J\lambda$, then
$
 \overline{v}(\mu) = \overline{v}_J(\mu)\overline{v}(\lambda)
$
 where the lengths add up. In particular, $\overline{v}(\lambda)^{-1}\in W_0^J$.
\end{Lemma}
\begin{proof}
 Let $\lambda_+\in L_+$ be the unique element in whose $W_0$-orbit $\lambda,\mu$ lie.

``$\le$'': We have
$
 \overline{v}_J(\mu)\overline{v}(\lambda)\lambda_+
 =\overline{v}_J(\mu)\lambda = \mu$,
 so we can write $\overline{v}_J(\mu)\overline{v}(\lambda)$ as $\overline{v}(\mu)w$ where $w$ fixes~$\lambda_+$. Since
 the stabiliser of $\lambda_+$ is a parabolic subgroup of $W_0$, and since $\overline{v}(\mu)$ has minimal length, in the above
 decomposition the lengths add up. Therefore, $\overline{v}(\mu)\le \overline{v}_J(\mu)\overline{v}(\lambda)$.

``$\ge$'':
 Let $v\in W_0^J, w\in W_J$ such that $\overline{v}(\mu)^{-1}=vw$,
 i.e., $w^{-1}v^{-1} = \overline{v}(\mu)$. Note that
 $v^{-1}\lambda_+\in L_{+,J}$ by Lemma~\ref{sec:lem-obtain-J-dominant}.
 Since $w^{-1}v^{-1}\lambda_+=\mu$, the element $v^{-1}\lambda_+$ lies in $W_J\mu$. Since it is $J$-dominant,
 it equals $\lambda$ by Proposition~\ref{sec:prop-J-fundamental-domain}.

 Thus, $w^{-1}\in W_J$ maps $\lambda$ to $\mu$, hence
 $\overline{v}_J(\mu)\le w^{-1}$ (same argument as for ``$\le$'').
 Furthermore, $v^{-1}$ maps $\lambda_+$ to $\lambda$, whence $\overline{v}(\lambda)\le v^{-1}$.
 Since the lengths of $v^{-1}$ and
 $w^{-1}$ add up, we can conclude
 \[
 \overline{v}_J(\mu)\overline{v}(\lambda)
 \le w^{-1}v^{-1}=\overline{v}(\mu),
 \]
 and hence the desired equality.

 To show that $v^{-1}=\overline{v}_J(\mu)$, note that modulo $W_J$ we have
$
 vW_J = \overline{v}(\lambda)^{-1}W_J$.
 Since $v$ is the shortest representative of this coset, we have $v\le \overline{v}(\lambda)^{-1}$.
 Together with $v^{-1}\ge \overline{v}(\lambda)$, we obtain $v=\overline{v}(\lambda)^{-1}\in W_0^J$.
\end{proof}

In order to later determine the leading coefficients of a $W_J$-invariant
polynomial, we now need to find out
which elements of a given $W_J$-orbit are the highest and lowest.
\begin{Corollary}
 Let $\lambda\in L_{+,J}$, then $\lambda$ is the smallest element of its $W_J$-orbit and
 $w_J\lambda$ is the largest.
\end{Corollary}
\begin{proof}
 As in Corollary~\ref{sec:cor-decomposition-symmetriser}, write
 $W_J^\lambda$ for the set of shortest representatives
 of~$W_{I,\lambda}$-cosets of~$W_J$ and let $w_{J,\lambda}$ be the longest element of~$W_{J,\lambda}$. As in the proof of Lemma~\ref{sec:lem-inverse-labelling}, the lengths of~$w_Jw_{J,\lambda}$ and
 $w_{J,\lambda}$ add up so that $w_Jw_{J,\lambda}$ is the shortest representative of $w_JW_{J,\lambda}$. Consequently, $\overline{v}_J(w_J\lambda) = w_Jw_{J,\lambda}$. Let
 $v\in W_J^\lambda$, then $v\le w_J$, so we can choose reduced expressions for $w_Jw_{J,\lambda}$ and $w_{J,\lambda}$ and
 obtain a reduced expression for $v$ by deleting reflections. Since $v$ is the shortest element in its
 right $W_{J,\lambda}$-coset, we need to delete all reflections belonging to $w_{J,\lambda}$, so that
 we have $v\le w_J w_{J,\lambda}$.

 Let now $\mu\in W_J\lambda$, then $\overline{v}_J(\mu)\in W_J^\lambda$ because it is the shortest element in its
 right $W_{J,\lambda}$-coset. Then we have $1\le \overline{v}_J(\mu)\le w_Jw_{J,\lambda}$, i.e.,
$
 \overline{v}_J(\lambda)\le \overline{v}_J(\mu) \le \overline{v}_J(w_J\lambda)$.
 We can multiply everything by $\overline{v}(\lambda)\in (W_0^J)^{-1}$ on the right since the lengths will always add up and
 obtain~${
 \overline{v}(\lambda)\le\overline{v}(\mu)\le\overline{v}(w_J\lambda)}
$
 by Lemma~\ref{sec:lem-decompose-vs}. By definition of the partial ordering, this implies that~${\lambda\le\mu\le w_J\lambda}$.
\end{proof}

\section{Intermediate Macdonald polynomials}\label{sec:poly}
We now begin by defining a class of polynomials possessing the desired
transformation behaviour under $\mathfrak{H}_J$.

\subsection{General}
\begin{Definition}
 For $\lambda\in L$, define \smash{$F^{(\epsilon)}_{J,\lambda} := U^{(\epsilon)}_J E_\lambda$}, where $E_\lambda$ is the non-symmetric Macdonald polynomial from
 Theorem~\ref{sec:thm-E}.
 Later, these will be suitably normalised to yield the intermediate Macdonald polynomials.
\end{Definition}

\begin{Corollary}
 Let $\lambda\in L$, $i\in J$, then
\smash{$ T_i F^{(\epsilon)}_{J,\lambda} = \tau_i^{(\epsilon)} F^{(\epsilon)}_{J,\lambda}$}.
 In particular, if $\epsilon$ is the trivial character, the polynomial
 \smash{$F^{(\epsilon)}_{J,\lambda}=F_{J,\lambda}$} is $W_J$-invariant.
\end{Corollary}
\begin{proof}
 Follows from applying Corollary~\ref{sec:cor-makes-symmetric} to
 the basic representation. If $\epsilon$ is trivial, we have%
 \[
 0 = (T_i-\tau_i)F^{(\epsilon)}_{J,\lambda}
 = (\bm{b}_{a_i} - \tau_i)(1-s_i)F^{(\epsilon)}_{J,\lambda}
 \]
 for all $i\in J$, and thus that \smash{$F^{(\epsilon)}_{J,\lambda}$} is $W_J$-symmetric.
\end{proof}

\begin{Corollary}
 Let $\lambda\in L$, $f\in A'_J$, then
 \[
 f(Y)F^{(\epsilon)}_{J,\lambda} = f(-r_{k'}(\lambda)) F^{(\epsilon)}_{J,\lambda},
 \]
 i.e., the commutative algebra $A'_J(Y)$ of difference-reflection operators
 diagonalises \smash{$F^{(\epsilon)}_{J,\lambda}$}.
\end{Corollary}
\begin{proof}
 $f(Y)$ commutes with $\mathfrak{H}_J$, hence also with \smash{$U^{(\epsilon)}_J$}. Then the claim follows from Theo\-rem~\ref{sec:thm-E}\,(ii).
\end{proof}

In some cases, the \smash{$F^{(\epsilon)}_{J,\lambda}$} turns out to be zero
(``when symmetrising something anti-symmetric''). We shall later see
that the following lemma indeed covers all instances where this happens.

\begin{Lemma}\label{sec:lem-F-vanishes}
 If there is $j\in J$ with $\epsilon(s_j)=-1$ and $s_j\lambda=\lambda$, then
 \smash{$F^{(\epsilon)}_{J,\lambda}=0$}.
\end{Lemma}
\begin{proof}
 Analogously to \cite[Section~5.7.1]{macdonald},
 using Lemma~\ref{sec:lem-symmetriser-symmetric}.
\end{proof}

Next, we will see that up to scalar multiples, \smash{$F^{(\epsilon)}_{J,\lambda}$}
only depends on
$\lambda$'s $W_J$-orbit. For that we define the labelling
$
 \epsilon k\colon S_J\to K, a\mapsto \epsilon(s_a)k(a)
$
of the root system $S_J$ spanned by $(a_j)_{j\in J}$, analogously $\epsilon k'$.

\begin{Lemma}\label{sec:lem-relation-within-orbit}
 Let $j\in J$ with $\langle \lambda, \alpha_j'\rangle>0$, then
 \[
 F^{(\epsilon)}_{J,s_j\lambda}
 = \epsilon(s_j) \tau_j \bm{c}'_{a'_j}(\epsilon(s_j)r_{k'}(\lambda))
 F^{(\epsilon)}_{J,\lambda}
 = \epsilon(s_j) \tau_j \bm{c}_{a'_j,\epsilon k'}(r_{k'}(\lambda))
 F^{(\epsilon)}_{J,\lambda}.
 \]
\end{Lemma}
\begin{proof}
 Analogously to \cite[Section~5.7.2]{macdonald}.
\end{proof}

This statement can be iterated.

\begin{Lemma}\label{sec:lem-relation-within-orbit-iterated}
 Let $\lambda\in L_{+,J}$ and $\mu\in W_J\lambda$, then
 \[
 F^{(\epsilon)}_{J,\mu}
 = \epsilon(\overline{v}_J(\mu))
 \tau_{\overline{v}_J(\mu)}
 \bm{c}_{\epsilon k', S'}(\overline{v}_J(\mu))(r_{k'}(\lambda))
 F^{(\epsilon)}_{J,\lambda},
 \]
 where as in {\rm \cite[equation~(4.4.6)]{macdonald}} we have
\smash{$
 \bm{c}_{k,S}(w) := \prod_{a\in S_1(w)}\bm{c}_{a,k}
$}
 for $w\in W$ $($analogously, for $k',S',w\in W')$.
\end{Lemma}
\begin{proof}
 Write $\overline{v}_J(\mu)=s_{i_1}\cdots s_{i_p}$ and define
$
 \mu_r := s_{i_{r+1}}\cdots s_{i_p} \lambda$,
 as well as $b'_r := s_{i_p}\cdots s_{i_{r+1}} a'_{i_r}$. By Lemma~\ref{sec:lem-flipped-roots}, we have
$
 S'_1(\overline{v}_I(\mu)) = \set{b'_1,\dots,b'_p}$.
 From Lemma~\ref{sec:lem-order-interval}, we furthermore know that
$
 \mu = \mu_0 > \cdots > \mu_p = \lambda$,
 so in particular the $\mu_r$ are all different and we have
$
 r_{k'}(\mu_r) = s_{i_{r+1}}\cdots s_{i_p} r_{k'}(\lambda)
$
 by iterating Lemma~\ref{sec:lem-pull-out-of-r}.
 This then shows that
 \begin{align*}
 \bm{c}_{a'_{i_r},\epsilon k'}(r_{k'}(\lambda_r))
 &= \bm{c}_{a'_{i_r}, \epsilon k'}(s_{i_{r+1}}\cdots s_{i_p} r_{k'}(\lambda))= \bm{c}_{s_{i_p}\cdots s_{i_{r+1}} a'_{i_r}, \epsilon k'}(r_{k'}(\lambda))= \bm{c}_{b'_r, \epsilon k'}(r_{k'}(\lambda)).
 \end{align*}
 We can now apply Lemma~\ref{sec:lem-relation-within-orbit} recursively
 to find
 \begin{align*}
 F^{(\epsilon)}_{J,\mu}
 &= \prod_{r=1}^p (\epsilon(s_{i_r}) \tau_{i_r} \bm{c}_{b'_r,\epsilon k'} (r_{k'}(\lambda)))
 F^{(\epsilon)}_{J,\lambda}= \epsilon(\overline{v}_J(\mu)) \tau_{\overline{v}_J(\mu)}
 \bm{c}_{\epsilon k', S'}(\overline{v}_J(\mu))(r_{k'}(\lambda))
 F^{(\epsilon)}_{J,\lambda}.\tag*{\qed}
 \end{align*}\renewcommand{\qed}{}
\end{proof}

\subsection{Leading terms}
Using results from the last section, we can now prove that the
\smash{$F^{(\epsilon)}_{J,\lambda}$} have non-vanishing leading coefficients unless
the conditions of Lemma~\ref{sec:lem-F-vanishes} are met.

\begin{Lemma}\label{sec:lem-F-leading}
 Let $\lambda\in L_{+,J}$ be such that $\epsilon(W_{J,\lambda})=1$. Then
 \[
 F^{(\epsilon)}_{J,\lambda} = \frac{W_{J,\lambda}\bigl(\tau_k^2\bigr)}{\tau_{w_J}} e(w_J\lambda) + \text{l.o.t.}
 \]
 In particular, $(F_{J,\lambda})_{\lambda\in L_{+,J}}$ is a basis of $A_J$.
\end{Lemma}
\begin{proof}
 By Lemma~\ref{sec:lem-Ti-action-E}\,(i), we have
 $T(w)E_\lambda=\tau_w E_\lambda$ for $w\in W_{J,\lambda}$, which shows that
 \[
 U_{J'} E_\lambda = \frac{W_{J,\lambda}\bigl(\tau_k^2\bigr)}{\tau_{w_{J,\lambda}}}
 E_\lambda,
 \]
 where $J'\subset J$ is the set of $i\in J$ such that $s_i\lambda=\lambda$.
 Using the decomposition from Corollary~\ref{sec:cor-decomposition-symmetriser}, we find
 \[
 F^{(\epsilon)}_{J,\lambda}
 = U^{(\epsilon)}_J E_\lambda
 = \frac{W_{J,\lambda}\bigl(\tau_k^2\bigr)}{\tau_{w_{J,\lambda}} \tau^{(\epsilon)}_{w_Jw_{J,\lambda}}}
 \sum_{v\in W_J^\lambda} \tau^{(\epsilon)}_v T(v) E_\lambda.
 \]
 From Lemma~\ref{sec:lem-Ti-action-E} and Theorem~\ref{sec:thm-E}, we
 conclude recursively that for all $v\in W_J^\lambda$ we have
 \[
 T(v)E_\lambda = \sum_{w\le v} c_{v,w} e(w\lambda)E_{w\lambda}
 = \sum_{w\le v}\sum_{\mu\le w\lambda} c_{v,w} d_{w\lambda,\mu} e(\mu)
 \]
 with $c_{v,v}=\tau_v^{-1}$ and $d_{\mu,\mu}=1$. Since for $w\in W_J$ we have
 $\overline{v}_J(w\lambda)\le w$, we have $w\lambda\le v\lambda$,
 and hence $\mu\le w\lambda\le v\lambda$ for every term $\mu$, $w$ occurring
 in this sum. Consequently, the highest term in $T(v)E_\lambda$ is
 $\tau_v^{-1}e(v\lambda)$, so that the leading coefficient of \smash{$F^{(\epsilon)}_{J,\lambda}$} comes only from $v=w_Jw_{J,\lambda}$, and equals
 \[
 \frac{W_{J,\lambda}\bigl(\tau_k^2\bigr)}{\tau_{w_{J,\lambda}} \tau^{(\epsilon)}_{w_Jw_{J,\lambda}}}
 \tau^{(\epsilon)}_{w_Jw_{J,\lambda}} \tau_{w_J w_{J,\lambda}}^{-1}
 = \frac{W_{J,\lambda}\bigl(\tau_k^2\bigr)}{\tau_{w_J}},
 \]
 which is nonzero.

 To see that the $F^{(\epsilon)}_{J,\lambda}$ form a basis, note that for $\lambda\in L_{+,J}$ we can define
 \[
 m_{J,\lambda}:= \sum_{\mu\in W_J\lambda} e(\mu),
 \]
 which is evidently a basis of $A_J$. In this basis, we have
 \[
 F_{J,\lambda} = \sum_{\substack{\mu\in L_{+,J}\\\mu\le\lambda}} c_{\lambda,\mu} m_{J,\mu}
 \]
 (note that for $\mu\le\lambda$ a nonzero term occurring in $F_{J,\lambda}$, the $J$-dominant
 element of $\mu$'s $W_J$-orbit~is lower than $\mu$, hence also $\le \lambda$) with
 \smash{$c_{\lambda,\lambda}=\frac{1}{\tau_{w_J}}W_{J,\lambda}\bigl(\tau_k^2\bigr)$}. This gives rise to a triangular matrix
 with invertible diagonals, so we can invert it (locally) and therefore find relations for the
 $m_{J,\lambda}$ in terms of $(F_{J,\mu})_{\mu\le\lambda}$.
\end{proof}

With this in hand we are ready to define the protagonist of this work.
\begin{Definition}\label{sec:def-P}
 Let
 \[
 P^{(\epsilon)}_{J,\lambda} := \frac{\tau_{w_J}}{W_{J,\lambda}\bigl(\tau_k^2\bigr)} F^{(\epsilon)}_{J,\lambda}.
 \]
 These are the \emph{intermediate Macdonald polynomials} associated to the parabolic subgroup
 \linebreak ${W_J \le W_0}$ and its character $\epsilon$.
\end{Definition}

By Lemma~\ref{sec:lem-F-leading} and Definition~\ref{sec:def-P}, we obtain
the following leading term.

\begin{Corollary}\label{sec:cor-leading-coeffs}
 If $\epsilon(W_{J,\lambda})=1$, we have
$\smash{P^{(\epsilon)}_{J,\lambda} = e(w_J \lambda)} +\text{l.o.t.}$
 In particular,
$
 P_{J,\lambda} = m_{J,\lambda} +\text{l.o.t.}$
 Otherwise, we have
\smash{$
 P^{(\epsilon)}_{J,\lambda} = 0$}.
\end{Corollary}

\subsection{Orthogonality}
As a consequence of the definition of the polynomials \smash{$F^{(\epsilon)}_{J,\lambda}$} in terms of symmetrisers, we have easy access to their
$A'_J(Y)$-eigenvalues. This allows us to prove their orthogonality.

\begin{Theorem}\label{sec:thm-ortho}
 For any $\epsilon$, the polynomials \smash{$\bigl(P^{(\epsilon)}_{J,\lambda}\bigr)_{\lambda\in L_{+,J}}$} from Definition~{\rm\ref{sec:def-P}} form a family of
 orthogonal polynomials with respect to the inner product from Definition~{\rm\ref{def-inner-product}}, with the caveat that \smash{$P^{(\epsilon)}_{J,\lambda}=0$} in
 case $-1\in\epsilon(W_{0,\lambda})$.
 In particular, the nonzero elements are an orthogonal basis of their
 span. For $\epsilon=1$, we thus obtain an orthogonal basis of $A_J$.
\end{Theorem}
\begin{proof}
 Let $\mu\ne \lambda\in L_{+,J}$. We will show that the corresponding inner product of $F$'s is zero for
 all choices of $k$ where all $k'\bigl(\alpha_i^\vee\bigr)>0$ for all $i\in I$.

 The weights $\lambda$, $\mu$ lie on different $W_J$-orbits. We first show that their $r_{k'}$s also lie on different $W_J$-orbits.
 Assume that there is $w\in W_J$ such that
 $wr_{k'}(\lambda)=r_{k'}(\mu)$. By definition, this means that
 \[
 wv(\lambda)^{-1}(\lambda_- - \rho_{k'})
 = v(\mu)^{-1}(\mu_- - \rho_{k'}).
 \]
 Since $\rho_{k'}$ is strictly dominant by Lemma~\ref{sec:lem-rho-inner-prod}, the vectors $\lambda_--\rho_{k'}$, $\mu_--\rho_{k'}$ are both
 strictly antidominant. Since they lie in the same $W_0$-orbit, they must be equal, hence
 $\lambda_-=\mu_-$ and hence $\mu$, $\lambda$ lie in the same $W_0$-orbit.
 Furthermore, their being strictly antidominant implies that they are
 regular, whence $wv(\lambda)^{-1} = v(\mu)^{-1}$. Consequently, we have
 \[
 w\lambda = wv(\lambda)^{-1}\lambda_-
 = v(\mu)^{-1}\mu_- = \mu,
 \]
 so that $\mu$, $\lambda$ are $W_J$-related. By Proposition~\ref{sec:prop-J-fundamental-domain}, we thus have $\mu=\lambda$, which is a contradiction.
 Therefore, $r_{k'}(\lambda)$, $r_{k'}(\mu)$ lie on different $W_J$-orbits.
 Hence the same is true for $-r_{k'}(\lambda)$, $-r_{k'}(\mu)$.

 Since $A'_J$ separates $W_J$-orbits of $V$, there exists $f\in A'_J$ such that $f(-r_{k'}(\lambda))\ne f(-r_{k'}(\mu))$. Then
 by~\eqref{eq:y-bernstein}, we have \smash{$f(Y)U^{(\epsilon)}_J=U^{(\epsilon)}_J$},
 so that
 \begin{align*}
 f(-r_{k'}(\lambda)) \bigl(F^{(\epsilon)}_{J,\lambda}, F^{(\epsilon)}_{J,\mu}\bigr) &=
 \bigl(f(-r_{k'}(\lambda)) F^{(\epsilon)}_{J,\lambda}, F^{(\epsilon)}_{J,\mu}\bigr)= \bigl(f(Y)F^{(\epsilon)}_{J,\lambda}, F^{(\epsilon)}_{J,\mu}\bigr)= \bigl(F^{(\epsilon)}_{J,\lambda}, f^*(Y)F^{(\epsilon)}_{J,\mu}\bigr)\\
 &= \bigl(F^{(\epsilon)}_{J,\lambda}, f^*(-r_{k'}(\mu)) F^{(\epsilon)}_{J,\mu}\bigr)= f(-r_{k'}(\mu)) \bigl(F^{(\epsilon)}_{J,\lambda}, F^{(\epsilon)}_{J,\mu}\bigr),
 \end{align*}
 whence \smash{$\bigl(F^{(\epsilon)}_{J,\lambda},F^{(\epsilon)}_{J,\mu}\bigr)=0$}.
 Consequently, we also have
$
 \bigl(F^{(\epsilon)}_{J,\lambda},F^{(\epsilon)}_{J,\mu}\bigr)=0$.

 By definition, the element \smash{$\bigl(P^{(\epsilon)}_{J,\lambda},P^{(\epsilon)}_{J,\mu}\bigr)_1$} is a rational function in finitely many values of $q$, and since the space of these values for positive $k$ is Zariski-dense, the equality \smash{$\bigl(P^{(\epsilon)}_{J,\lambda},P^{(\epsilon)}_{J,\mu}\bigr)=0$} holds identically.
\end{proof}

\begin{Example}
 For $\epsilon=1$ and $J=I_0$, the leading coefficient from
 Corollary~\ref{sec:cor-leading-coeffs} together with the orthogonality
 from Theorem~\ref{sec:thm-ortho} implies that
 \begin{itemize}\itemsep=0pt
 \item $\forall\mu\in L_+\colon P_{J,\mu} = m_\mu + \text{l.o.t.}$
 \item $\forall\lambda,\mu\in L_+\colon \lambda<\mu\Rightarrow (P_{J,\mu}, m_\lambda)=0$.
 \end{itemize}
 Together with \cite[Sections 5.1.35 and 5.3.1f]{macdonald}, this implies
 that our intermediate Macdonald polynomials are just the symmetric
 Macdonald polynomials from \cite[Section~5.3]{macdonald}.

 For $\epsilon=1$ and $J=\varnothing$, our definitions show that
 $P_{J,\lambda}=E_\lambda$, i.e., that the intermediate Macdonald polynomials
 are the nonsymmetric Macdonald polynomials.

 These two examples illustrate the choice of the name ``intermediate Macdonald polynomials''.
\end{Example}

\subsection{Norms}
Using an approach similar to \cite[Section~5.7.12]{macdonald}, we can now also
express the norm of \smash{$P^{(\epsilon)}_{J,\lambda}$} in terms of
the norm of $E_{w_J\lambda}$.

\begin{Theorem}
 If $\epsilon(W_{J,\lambda})=1$, we have
 \[
 \bigl(P^{(\epsilon)}_{J,\lambda}, P^{(\epsilon)}_{J,\lambda}\bigr)
 = \epsilon(w_J) \frac{\bm{c}_{S', -\epsilon k'}(w_Jw_{J,\lambda})(r_{k'}(\lambda)) W_J^\lambda\bigl(\tau^{(\epsilon)2}_k\bigr)}{\tau_{w_Jw_{J,\lambda}}^2\tau^{(\epsilon)}_{w_Jw_{J,\lambda}}}
 (E_{w_J\lambda}, E_{w_J\lambda}).
 \]
\end{Theorem}
\begin{proof}
 We have
 \begin{align*}
 \bigl(P^{(\epsilon)}_{J,\lambda}, P^{(\epsilon)}_{J,\lambda}\bigr)
 = \frac{\tau_{w_{J,\lambda}}^2}{W_{J,\lambda}\bigl(\tau_k^2\bigr)^2}
 \bigl(F^{(\epsilon)}_{J,\lambda}, F^{(\epsilon)}_{J,\lambda}\bigr).
 \end{align*}
 By Lemma~\ref{sec:lem-relation-within-orbit-iterated}, this equals
 \[
 \epsilon(w_Jw_{J,\lambda}) \frac{\tau_{w_{J,\lambda}}^2 \tau_{w_Jw_{J,\lambda}}}{W_{J,\lambda}\bigl(\tau_k^2\bigr)^2\bm{c}_{-\epsilon k',S'}(w_Jw_{J,\lambda})(r_{k'}(\lambda))}
 \bigl(F^{(\epsilon)}_{J,\lambda}, F^{(\epsilon)}_{J,w_J\lambda}\bigr).
 \]
 Note that $\epsilon(w_{J,\lambda})=1$, so that $\epsilon(w_Jw_{J,\lambda})=\epsilon(w_J)$.
 Applying Proposition~\ref{sec:prop-symmetriser-properties}\,(iii), we get
 \[
 = \epsilon(w_J)
\frac{W_J\bigl(\tau^{(\epsilon)2}_k\bigr)\tau_{w_{J,\lambda}}^2 \tau_{w_Jw_{J,\lambda}}}{W_{J,\lambda}\bigl(\tau_k^2\bigr)^2\tau^{(\epsilon)}_{w_J}\bm{c}_{-\epsilon k',S'}(w_Jw_{J,\lambda})(r_{k'}(\lambda))}
 \bigl(F^{(\epsilon)}_{J,\lambda}, E_{w_J\lambda}\bigr).
 \]
 Since the $e(w_J\lambda)$-coefficient (the leading coefficient)
 in $F^{(\epsilon)}_{J,\lambda}$ is $\frac{W_{J,\lambda}(\tau_k^2)}{\tau_{w_J}}$,
 we have
 \[
 F^{(\epsilon)}_{J,\lambda} = \sum_{\mu\in W_J\lambda} c_\mu E_\mu
 \]
 with \smash{$c_{w_J\lambda} = \frac{W_{J,\lambda}(\tau_k^2)}{\tau_{w_J}}$}. Since all
 other $E_\mu$ are orthogonal to $E_{w_J\lambda}$, we have
 \[
 \bigl(F^{(\epsilon)}_{J,\lambda}, E_{w_J\lambda}\bigr)
 = \frac{W_{J,\lambda}\bigl(\tau_k^2\bigr)}{\tau_{w_J}} (E_{w_J\lambda}, E_{w_J\lambda}),
 \]
 so that
 \[
 \bigl(P^{(\epsilon)}_{J,\lambda}, P^{(\epsilon)}_{J,\lambda}\bigr)
 = \epsilon(w_J) \frac{W_J\bigl(\tau^{(\epsilon)2}_k\bigr)\tau_{w_{J,\lambda}} \tau_{w_Jw_{J,\lambda}}}{W_{J,\lambda}\bigl(\tau_k^2\bigr)\tau^{(\epsilon)}_{w_J}\bm{c}_{-\epsilon k',S'}(w_Jw_{J,\lambda})(r_{k'}(\lambda))}
 (E_{w_J\lambda}, E_{w_J\lambda}).
 \]
 Next, note that $\tau_{w_{J,\lambda}}=\tau^{(\epsilon)}_{w_{J,\lambda}}$,
 and that \smash{$W_{J,\lambda}\bigl(\tau_k^2\bigr)=W_{J,\lambda}\bigl(\tau^{(\epsilon)2}_k\bigr)$},
 so that
 \[
 \bigl(P^{(\epsilon)}_{J,\lambda}, P^{(\epsilon)}_{J,\lambda}\bigr)
 = \epsilon(w_J) \frac{W_J^\lambda\bigl(\tau^{(\epsilon)2}_k\bigr)\tau_{w_J w_{J,\lambda}}}{\tau^{(\epsilon)}_{w_Jw_{J,\lambda}}
 \bm{c}_{-\epsilon k', S'}(w_Jw_{J,\lambda})(r_{k'}(\lambda))}
 (E_{w_J\lambda}, E_{w_J\lambda}).\tag*{\qed}
 \] \renewcommand{\qed}{}
\end{proof}

\section{Invariant vector-valued polynomials}\label{sec:jasper}
A method of obtaining a representation of $\tilde{\mathfrak{H}}$ on
vector-valued polynomials is presented in~\cite{cherednikInduced}. Given a
$\mathfrak{H}$-module $V$, we can consider the $\tilde{\mathfrak{H}}$-module
$\tilde{\mathfrak{H}}\otimes_{\mathfrak{H}}V$. But which module $V$ should
we start with? Taking inspiration from \cite[Section~3]{stokmanInduced},
we consider principal series modules.

\subsection[Y-parabolically induced modules of H]{$\boldsymbol{ Y}$-parabolically induced modules of $\boldsymbol{\mathfrak{H}}$}
\label{sec:parabolically-induced}
\begin{Lemma}
 The set $\mathfrak{H}_J A'(Y)$ is a subalgebra of $\mathfrak{H}$.
\end{Lemma}
\begin{proof}
 Obviously $\mathfrak{H}_J$ and $A'(Y)$ are subalgebras of their own, and
 the Bernstein--Lusztig--Zelevinsky presentation~\eqref{eq:y-bernstein}
guarantees that their product forms an algebra.
\end{proof}

\begin{Definition}
 Let $\phi\colon \mathfrak{H}_J A'(Y)\to K$ be a ring homomorphism.
 We can define the following induced modules with it:
 \begin{enumerate}\itemsep=0pt
 \item[(i)] The \emph{principal series representation} is
 the $\mathfrak{H}$-module
$
 M_J(\phi) \!\!:=\! \mathfrak{H} \otimes_{\mathfrak{H}_JA'(Y)}\! K
$
(where $\mathfrak{H}_JA'(Y)$ acts on $K$ via $\phi$);
\item[(ii)] the \emph{standard $Y$-parabolically induced $\tilde{\mathfrak{H}}$-module} is the $\tilde{\mathfrak{H}}$-module
 \[
 \mathbb{M}_J(\phi) := \tilde{\mathfrak{H}}\otimes_{\mathfrak{H}}
 M_J(\phi) = \tilde{\mathfrak{H}}\otimes_{\mathfrak{H}_JA'(Y)} K.
 \]
 \end{enumerate}
 Moreover, an element $v\in\mathbb{M}_J(\phi)$ is called \emph{spherical} if
 $(T_i-\tau_i) v = 0$ for all $i\in I_0$. Write~$\mathbb{M}_J(\phi)^{\mathfrak{H}_0}$ for the vector space of spherical
 vectors.
\end{Definition}

\begin{Lemma}\label{sec:lem-induced-module-basis}
 A $K$-basis of $\mathbb{M}_J(\phi)$ is given by
$
 (X^\mu T(v))_{\mu\in L, v\in W_0^J}$.
 In particular, $\mathbb{M}_J(\phi)$ is a~free $A$-module
 with basis $(T(v))_{v\in W_0^J}$.
\end{Lemma}
\begin{proof}
 By Lemma~\ref{sec:lem-basis-daha}, we can decompose
$\tilde{\mathfrak{H}} \cong A(X)\otimes\mathfrak{H}_0\otimes A'(Y)$
 as right $A'(Y)$-modules. Note that using the shortest coset representatives,
 we can also decompose
 \[
 \mathfrak{H}_0 = \bigoplus_{v\in W_0^J} T(v) \mathfrak{H}_J
 \]
 as right $\mathfrak{H}_J$-modules. Consequently,
 \[
 \tilde{\mathfrak{H}} \cong A(X)\otimes \bigoplus_{v\in W_0^J} KT(v)\otimes \mathfrak{H}_J A'(Y)
 \]
 as right $\mathfrak{H}_J A'(Y)$-modules. Consequently,
 the family \smash{$(X^\mu T(v))_{\mu\in L, v\in W_0^J}$} is a basis of the
 right $\mathfrak{H}_J A'(Y)$-module $\tilde{\mathfrak{K}}$, and hence a
 $K$-basis of $\mathbb{M}_J(\phi)$.
\end{proof}

Purely on a vector space level, we can therefore establish
a bijection $A\otimes K[W_0/W_J]\cong \mathbb{M}_J(\phi)$ by mapping
$f\otimes \delta_{vW_J}\mapsto f(X)T(v)$ where $v\in W_0^J$.

\begin{Example}
 For $J=I_0$, the $Y$-parabolic subalgebra $\mathfrak{H}_JA'(Y)$ is just
 $\mathfrak{H}$. For the map $\phi$ from
 Lemma~\ref{sec:lem-trivial-character}, the standard $Y$-parabolically
 induced $\tilde{\mathfrak{H}}$-module $\mathbb{M}_J(\phi)$ is how we
 defined the basic representation on~$A$.
\end{Example}

\subsection{Spherical vectors}\label{sec:spherical-vec}
We shall assume now that our character $\phi$ restricts to the
trivial character of $\mathfrak{H}_J$. Then the induced module
$\mathbb{M}_J(\phi)$ corresponds to $(\C[P]\otimes\C[W/W_J])^W$ from
\cite{mvpVectorValued} in the classical limit. We will see that the result \cite[Lemma~5.1]{mvpVectorValued} carries over to the $q$-setting and the spherical vectors
of~$\mathbb{M}_J(\phi)$
correspond to elements of $A_J$.

\begin{Definition}
 We define $\Gamma\colon A_J\to \mathbb{M}_J(\phi)$ by mapping
$ f \mapsto U_0 f(X)$,
 where $U_0$ refers to the symmetriser for the whole group $W_0$.
\end{Definition}

\begin{Lemma}\label{sec:lem-Gamma-invariant-highest-term}
 The map $\Gamma$ is a well-defined $A_0$-linear map
 $A_J\to (\mathbb{M}_J(\phi))^{\mathfrak{H}_0}$.
 In particular, we have
 \[
 \Gamma(f) = \sum_{v\in W_0^J} f_v(X) T(v),\qquad f_{w_0w_J}=\frac{W_J\bigl(\tau_k^2\bigr)}{\tau_{w_J}} w_0f.
 \]
\end{Lemma}
\begin{proof}
 By Lemma~\ref{sec:lem-symmetriser-symmetric}, we have $(T_i-\tau_i)U_0=0$ for all $i\in I_0$. Consequently,
 by associativity also
$0 = (T_i-\tau_i)U_0 f(X) = (T_i-\tau_i)\Gamma(f)$.
 For $A_0$-linearity, note that $A_0(X)$ commutes with
 $\mathfrak{H}_0$, so that it also commutes with $U_0$.

 For the second claim, note that by Corollary~\ref{sec:lem-decomposition-symmetriser} we have
 \[
 \Gamma(f) = \tau_{w_0w_J}^{-1} \sum_{v\in W_0^J} \tau_v T(v)
 U_J f(X).
 \]
 Since $f\in A_J$, the elements $f(X)$ and $U_J$ commute, and we have
 \[
 U_J = \frac{W_J\bigl(\tau_k^2\bigr)}{\tau_{w_J}}
 \]
 in $\mathbb{M}_J(\phi)$. This shows that
 \begin{align*}
 \Gamma(f)
 &= \frac{W_J\bigl(\tau_k^2\bigr)}{\tau_{w_0}} \sum_{v\in W_0^J} \tau_v T(v)
 f(X)= \frac{W_J\bigl(\tau_k^2\bigr)}{\tau_{w_0}} \sum_{v\in W_0^J} \tau_v
 \sum_{w\le v} f_{v,w}(X) T(w),
 \end{align*}
 where $f_{v,v} = v f$ (follows inductively from~\eqref{eq:x-bernstein}). Since $w_0w_J$ is the smallest element of its right
 $W_J$-coset, the only summand contributing to the $T(w_0w_J)$ term
 is where $v=w=w_0w_J$. As a~consequence, when we expand $\Gamma(f)$ in
 terms of our basis from Lemma~\ref{sec:lem-induced-module-basis},
 say
 \[
 \Gamma(f) = \sum_{v\in W_0^J} f_v(X)T(v)
 \]
 with
 \[
 f_{w_0w_J} = \frac{W_J\bigl(\tau_k^2\bigr)}{\tau_{w_0}} \tau_{w_0w_J}
 w_0w_J f = \frac{W_J\bigl(\tau_k^2\bigr)}{\tau_{w_J}} w_0 f.\tag*{\qed}
 \]\renewcommand{\qed}{}
\end{proof}

\begin{Remark}\label{sec:rk-htilde-intertwiner}
 It is tempting to view $\Gamma$ as the concatenation of
 the $\tilde{\mathfrak{H}}$-linear map $\iota\colon A\to\mathbb{M}_J(\phi)$ mapping $1\mapsto 1$ and multiplication with $U_0$, rendering $\Gamma$ in fact
 not only $A(X)$-linear but indeed $Z_{\tilde{\mathfrak{H}}}(U_0)$-linear \big(the centraliser of $U_0$ in $\tilde{\mathfrak{H}}$\big).
 This, however, falls flat when we try to investigate if the map $\iota$
 even exists in the first place. Since we defined both \smash{$\tilde{\mathfrak{H}}$}-module structures
 in terms of induced representations, we can map $1\mapsto 1$ if
 \smash{$\operatorname{Ann}_A(1)\subset\operatorname{Ann}_{\mathbb{M}_J(\phi)}$}.
 Since $1$ lies in the centre of~$\tilde{\mathfrak{H}}$, the annihilators
 are the kernels of~$\tilde{\phi}$, $\phi$, respectively (here, we write
 $\tilde{\phi}$ for the trivial character~$\phi$ from Lemma~\ref{sec:lem-trivial-character}). However, unless
 $J=I_0$, the kernel of $\tilde{\phi}$ is bigger
 than that of~$\phi$.
\end{Remark}

\begin{Corollary}\label{sec:cor-Gamma-injective}
 $\Gamma$ is injective.
\end{Corollary}
\begin{proof}
 Since \smash{$\frac{W_J(\tau_k^2)}{\tau_{w_J}}$} is nonzero and
 $w_0$ is an automorphism of $A$, Lemma~\ref{sec:lem-Gamma-invariant-highest-term} implies that we can read off $f$ from $\Gamma(f)$'s
 expansion in the basis from Lemma~\ref{sec:lem-induced-module-basis}.
\end{proof}

If we want to show that $\Gamma$ is also surjective (onto the spherical
vectors of $\mathbb{M}_J(\phi)$), we need to find out more about the
spherical vectors, especially how the coefficients are related.

\begin{Proposition}\label{sec:prop-Gamma-properties}
 Let $h\in\mathbb{M}_J(\phi)$ be spherical, say
 \[
 h = \sum_{v\in W_0^J} f_v(X) T(v),
 \]
 then
 \begin{enumerate}\itemsep=0pt
 \item[$(i)$] Then $h$ is uniquely determined by any $f_v$.
 \item[$(ii)$] $f_{w_0w_J}$ is $w_0 W_J w_0$-symmetric.
 \item[$(iii)$] $w_0 f_{w_0 w_J}\in A_J$.
 \end{enumerate}
\end{Proposition}
\begin{proof}
 We begin by describing the action of $\mathfrak{H}_0$ on $M_J(\phi)$.
 If $v\in W_0^J$ and $i\in I_0$, we can decompose $s_iv\in W_0$ according
 to $W_0^J W_J$, say as
$s_i v =: (s_i\bullet v) m_i(v)$.
 Then we have $s_i vW_J= (s_i\bullet v)W_J$, so that $\bullet$ describes
 the $W_0$-action on $W_0/W_J$ on shortest coset representatives,
 and the $m$'s provide the cocycles. If $\ell(s_iv)>\ell(v)$, we have
 \[
 T_i T(v) = T(s_iv) = T((s_i\bullet v) m_i(v))
 = T(s_i\bullet v) T(m_i(v))
 = \tau_{m_i(v)} T(s_i\bullet v)
 \]
 in $M_J(\phi)$. Otherwise, we have
 \[
 T_i T(v) = T(s_iv) + \bigl(\tau_i - \tau_i^{-1}\bigr)
 T(v)
 = \tau_{m_i(v)} T(s_i\bullet v)
 + \bigl(\tau_i - \tau_i^{-1}\bigr) T(v).
 \]
 In the special case of $v=w_0w_J$, and $s_i\in w_0W_Jw_0$, we have
 $w_0 s_i w_0\in W_J$, and since the longest element of a Coxeter group
 always permutes the simple reflections (by conjugation), we have that
 $w_J w_0 s_i w_0 w_J = s_j$ is a simple reflection in $W_J$. This implies
$
 s_i w_0 w_J = w_0 w_J s_j$,
 which by Lemma~\ref{sec:permute-simple-reflections} implies that
 \[
 T_i T(w_0w_J) = T(w_0w_J) T_j = \tau_j T(w_0w_J)
 = \tau_i T(w_0w_J).
 \]

(i) We now consider $T_i$'s action on $h$. By the Bernstein--Lusztig--Zelevinsky
 presentation~\eqref{eq:x-bernstein} and the action of $T_i$ on $M_I(\phi)$, we obtain
 \begin{align*}
 T_i h ={}& \sum_{v\in W_0^J} (s_i f_v)(X) T_i T(v)
 + \bm{b}_i(X) \sum_{v\in W_0^J} (f_v-s_if_v)(X) T(v)\\
 ={}& \sum_{v\in W_0^J} (s_i f_v)(X) \tau_{m_i(v)} T(s_i\bullet v)
 + \bm{b}_i(X) \sum_{v\in W_0^J} (f_v-s_if_v)(X) T(v)\\
 &+ \bigl(\tau_i-\tau_i^{-1}\bigr)\sum_{\substack{v\in W_0^J\\\ell(s_iv)<\ell(v)}} f_v(X) T(v)\\
 ={}& \sum_{v\in W_0^J} (\tau_{m_i(s_i\bullet v)} s_i f_{s_i\bullet v} + \bm{b}_i (f_v-s_if_v))(X) T(v)\\
 &+ \bigl(\tau_i-\tau_i^{-1}\bigr)\sum_{\substack{v\in W_0^J\\\ell(s_iv)<\ell(v)}} (s_if_v)(X) T(v).
 \end{align*}
 Since the $(T(v))_{v\in W_0^J}$ are $A(X)$-linearly independent, we have
 \[
 \tau_i f_v = \tau_{m_i(s_i\bullet v)}
 s_i f_{s_i\bullet v} + \bm{b}_i(f_v-s_if_v)
 + \begin{cases}
 0, & \ell(s_iv)>\ell(v),\\
 \bigl(\tau_i-\tau_i^{-1}\bigr)(s_if_v), & \ell(s_iv)<\ell(v),
 \end{cases}
 \]
 We can solve for $f_{s_i\bullet v}$ and thus obtain that
 $f_{s_i\bullet v}$ is uniquely determined by $f_v$. Inductively, we thus
 obtain that $f_{w\bullet v}$ is uniquely determined by $f_v$ for all
 $w\in W_0$. Since $W_0$ acts transitively on $W_0/W_J$ and hence also
 on $W_0^J$, this shows that all coefficient polynomials~$(f_v)_{v\in W_0^J}$ are
 determined by any one $f_v$.

(ii) If we consider the case of $s_i\in w_0W_Jw_0$ and $v=w_0w_J$, we find
 that $s_i\bullet v=v$, $m_i(v)=s_i$, and $\ell(s_iv)>\ell(v)$, so that
$
 \tau_i f_v = \tau_i (s_if_v) + \bm{b}_i(f_v-s_if_v)$,
 in other words,
$
 0 = (\bm{b}_i-\tau_i)(f_v-s_if_v)$,
 which implies that $f_v=s_if_v$. In other words, $f_{w_0w_J}$ is
 $w_0W_Jw_0$-symmetric.

(iii) Let $w\in W_J$, then
$
 w w_0 f_{w_0w_J} = w_0 w_0 w w_0 f_{w_0w_J}
 = w_0 f_{w_0w_J}
$
 as $w_0ww_0\in w_0W_Jw_0$.
\end{proof}

\begin{Theorem}
 $\Gamma\colon A_J\to \mathbb{M}_J(\phi)^{\mathfrak{H}_0}$ is an
 $A_0$-linear isomorphism.
\end{Theorem}
\begin{proof}
 From Lemma~\ref{sec:lem-Gamma-invariant-highest-term}, we know that
 $\Gamma$ is well-defined and $A_0$-linear. With Corollary~\ref{sec:cor-Gamma-injective}, we conclude that $\Gamma$ is injective. For surjectivity,
 let
 \[
 h = \sum_{v\in W_0^J} f_v(X)T(v)\in\mathbb{M}_J(\phi)^{\mathfrak{H}_0}.
 \]
 From Proposition~\ref{sec:prop-Gamma-properties}\,(iii), we know that
 $w_0f_{w_0w_J}\in A_J$. Consider now
 \[
 h' := \Gamma\qty(\frac{\tau_{w_J}}{W_J\bigl(\tau_k^2\bigr)} w_0f_{w_0w_J})
 = \sum_{v\in W_0^J} f'_v(X)T(v).
 \]
 By Lemma~\ref{sec:lem-Gamma-invariant-highest-term}, we know that
 \[
 f'_{w_0w_J} = \frac{W_J\bigl(\tau_k^2\bigr)}{\tau_{w_J}}
 \frac{\tau_{w_J}}{W_J\bigl(\tau_k^2\bigr)} w_0 w_0 f_{w_0w_J}
 = f_{w_0w_J}.
 \]
 By Proposition~\ref{sec:prop-Gamma-properties}\,(i) therefore,
 all other coefficient polynomials of $h$, $h'$ also have to be equal,
 whence $h=h'$.
\end{proof}

Using this isomorphism, we can now push $(\cdot,\cdot)_1$ to
$\mathbb{M}_J(\phi)^{\mathfrak{H}_0}$ and obtain the basis of vector-valued
orthogonal polynomials $(\Gamma(P_{J,\lambda}))_{\lambda\in L_{+,J}}$.

In \cite[Section~6.6]{stokmanInduced}, the authors introduce a basis
of $W_0$-symmetric orthogonal elements of~$\mathbb{M}_J(\phi)$ that
diagonalise the action of $A'_0(Y)$. Since our map $\Gamma$ is not
expected to preserve
any of the $Y$-actions, we do not expect this basis to coincide with the intermediate Macdonald polynomials.
However, in \cite[Theorem~7.6]{Ven25} it is shown that they coincide in the
$q\to\infty$ limit.

\section{Vector-valued invariant polynomials}\label{sec:w0-invariant}
Classically, there exists an isomorphism
\[
 (A\otimes K[W_0/W_J])^{W_0}
 \cong A^{W_0} \otimes K[W_0/W_J]
 = A_0^{\# W_0^J}.
\]
The reason for this is that as is shown in~\cite{pittieSteinberg},
$A_J$ is a free $A_0$-module. This is proved by constructing an appropriate
basis $(e_v)_{v\in W_0^J}$, where $e_1=1$. Note that this
particular choice of basis is by no means unique or canonical, and we shouldn't
restrict ourselves to it.

\begin{Lemma}\label{sec:lem-matrix-weight}
 Let \smash{$(e_v)_{v\in W_0^J}$} be an $A_0$-basis of $A_J$. Via this basis, we identify
 elements $f\in A_J$ with column vectors $\underline{f}$ $\big($indexed by $W_0^J\big)$. Define
 \[
 m_{v,v'} := \frac{1}{\#W_0} \sum_{w\in W_0} w\frac{e_v e_{v'}^*}{\Delta_0} \in K(L)^{W_0}
 \]
 and write $M=(m_{v,v'})_{v,v'\in W_0^J}$ for the matrix thus defined. Then,
\smash{$
 (f,g) = \operatorname{ct}\qty(\underline{f}^{\mathsf T} M \underline{g}^* \nabla)
$}
 for all $f,g\in A_J$. Here, $\nabla$ is the $W_0$-invariant weight
 function defined in Definition~{\rm \ref{sec:def-deltas-nabla}}.
\end{Lemma}
\begin{proof}
 If we expand
 \[
 f = \sum_{v\in W_0^J} f_v e_v, \qquad g = \sum_{v\in W_0^J} g_v e_v,
 \]
 then the claim states that
 \[
 (f,g) = \sum_{v,v'\in W_0^J} \operatorname{ct}\qty(f_v m_{v,v'} g_{v'}^*\nabla).
 \]
 To prove this, we substitute the expansions of $f$, $g$
 \[
 (f,g) = \sum_{v,v'\in W_0^J} \operatorname{ct}\qty(f_v g_{v'}^*
 \nabla \frac{e_v e_{v'}^*}{\Delta_0}).
 \]
 Recall that we can write $\Delta=\nabla\Delta_0^{-1}$ where $\nabla$ is
 $W_0$-invariant by Lemma~\ref{lem-weight-distributions}\,(ii).
 Next, note that the notion of the constant term is $W_0$-invariant, as is
 $f_v g_{v'}^*\nabla$, so the summands are not changed by also symmetrising over
 $W_0$
 \[
 (f,g) = \sum_{v,v'\in W_0^J} \operatorname{ct}\qty(
 f_v g_{v'}^*\nabla \frac{1}{\#W_0} \sum_{w\in W_0} w \frac{e_v e_{v'}^*}{\Delta_0})
 = \sum_{v,v'\in W_0^J} \operatorname{ct}\qty(f_v m_{v,v'} g_{v'}^* \nabla).\tag*{\qed}
 \] \renewcommand{\qed}{}
\end{proof}

\begin{Corollary}
 For any fixed basis $(e_v)_{v\in W_0^J}$, the family
 $(\underline{P_{J,\lambda}})_{\lambda\in L_{+,J}}$ is a family of vector-valued polynomials
 that is orthogonal with respect to the matrix weight $\nabla M$.
\end{Corollary}

We close this discussion of general theory with some observations about
matrix weights.

\begin{Definition}
 A matrix weight $M\in (K(L)^{W_0})^{n\times n}$ is called \emph{reducible} if there is
 $R\in\operatorname{GL}_n(K)$ such that $R M R^{*{\mathsf T}}$ is a block matrix.
(Here, $*$ refers to entry-wise application of $K$'s $*$ involu\-tion.)
\end{Definition}

\begin{Corollary}\label{sec:cor-shi-block-structure}
 If $M$ is a matrix weight whose entries are polynomials,
 say $M= \sum_{\mu\in L_+} M_\mu m_\mu$, then $M$ is reducible if all
 $M_\mu$ can be brought into the same block shape by a similarity transformation
 with the same matrix~$R$.
 In particular, $M$ can be made diagonal iff all the $M_\mu$ can be made diagonal
 with the same matrix~$R$.
\end{Corollary}
\begin{proof}
 If $D=RMR^{*{\mathsf T}}$ is a block matrix, we can expand $D$ as a polynomial matrix as
 \[
 D = \sum_{\mu\in L_+} R M_\mu R^{*{\mathsf T}} m_\mu.
 \]
 Since the $m_\mu$ are $K$-linearly independent, all $R M_\mu R^{*{\mathsf T}}$ have to
 have the same block structure as~$D$.
\end{proof}

We now consider two examples, where we are able to
transform the basis constructed in~\cite{pittieSteinberg} by a
unipotent triangular matrix (over $A_0$) to obtain a maximally reduced
matrix weight.

\subsection[Example: (C\_1\^{}vee, C\_1) with J=varnothing]{Example: $\boldsymbol{\bigl(C_1^\vee, C_1\bigr)}$ with $\boldsymbol{J=\varnothing}$}
This example was taken (with slightly different conventions) from
\cite{koornwinderBouzeffour} and \cite[Section~6]{koornwinderMazzocco}. The~${q\to1}$ limit can be compared to~\cite{mvpBC1}. For
an affine root system of type $\bigl(C_1^\vee, C_1\bigr)$, the double-affine Hecke
algebra $\tilde{\mathfrak{H}}$ is
generated by $T_1$, $X$, $Y$ (in the notation of \cite[Section~6.4]{macdonald}) and depends on the four parameters
$k_1$, $k_2$, $k_3$, $k_4$. As is customary, we use these four parameters to define
the more convenient Askey--Wilson parameters
\[
 (a,b,c,d) = \bigl(q^{k_1},-q^{k_2},q^{1/2+k_3},-q^{1/2+k_4}\bigr).
\]
We pick $J=\varnothing$ and thus represent $A$ as a free $A_0$-module.
Steinberg's proof provides the basis $\tilde{e}_1=1$, $\tilde{e}_{s_1} = x^{-1}$ (writing $x:= e(a_1)$).
We shall first consider this basis and compare our findings to
\cite{mvpBC1}.

\begin{Lemma}\label{sec:lem-weight-matrix-steinberg}
 In the basis $(\tilde{e}_1,\tilde{e}_{s_1})$, we obtain $($as in Lemma~{\rm \ref{sec:lem-matrix-weight})}
 the following
 weight matrix:
 \[
 \tilde{M} = \frac{1}{2}\mqty(1-ab & x+x^{-1}-a-b\\-ab\bigl(x+x^{-1}\bigr)+a+b & 1-ab).
 \]
 We can now do a similarity transform with the constant matrix
 \[
 U:=\mqty(-a & 1\\-b & 1)
 \]
 to find
 \[
 U\tilde{M}U^{*{\mathsf T}} = \frac{a-b}{2}\mqty(-\bigl(x+x^{-1}\bigr)+a+a^{-1} & 0\\0 & x+x^{-1}-b-b^{-1}).
 \]
\end{Lemma}
\begin{proof}
 Note that
 \[
 \Delta_0 = \frac{x-x^{-1}}{(x-a)\bigl(1-bx^{-1}\bigr)}
 = \frac{\delta}{F},
 \]
 where $\delta=x-x^{-1}$ is the Weyl-denominator, and $F=x\bigl(1-ax^{-1}\bigr)\bigl(1-bx^{-1}\bigr)$. Then
 \[
 \tilde{m}_{v,v'} = \frac{1}{2\delta}(1-s_1)(F\tilde{e}_v\tilde{e}_{v'}^*)
 \]
 for all $v,v'\in W_0$
 since $s_1\delta^{-1} = -\delta^{-1}s_1$. This means that
 \begin{align*}
 \tilde{m}_{1,1}&=\tilde{m}_{s_1,s_1}=\frac{1}{2\delta}(1-s_1)F\\
 &= \frac{1}{2\delta}\qty(x - a - b + abx^{-1}
 - x + a + b - abx)= \frac{1-ab}{2},\\
 \tilde{m}_{1,s_1}&=\frac{1}{2\delta}(1-s_1)(Fx)= \frac{1}{2\delta}\qty(x^2 - ax-bx+ab - x^{-2} + ax^{-1}+bx^{-1}-ab)\\
 &= \frac{x+x^{-1}-a-b}{2},\\
 \tilde{m}_{s_1,1}&=\frac{1}{2\delta}(1-s_1)(Fx^{-1})= \frac{1}{2\delta}\qty(1 - ax^{-1}-bx^{-1} + abx^{-2}
 - 1 + ax + bx - abx^2)\\
 &= \frac{-ab\bigl(x+x^{-1}\bigr) + a + b}{2}.
 \end{align*}
 Lastly, to see the similarity transform, note that
 \[
 2\tilde{M} = \mqty(1-ab & -a-b\\a+b & 1-ab)
 + \mqty(0 & 1\\-ab & 0)\bigl(x+x^{-1}\bigr)
 = 2\tilde{M}_0 m_0 + 2\tilde{M}_{a_1} m_{\alpha_1}
 \]
 (using notation from Corollary~\ref{sec:cor-shi-block-structure}).
 We have
 \begin{align*}
 U2\tilde{M}_0U^{*{\mathsf T}} &= \mqty(-a & 1\\-b & 1)\mqty(1-ab & -a-b\\a+b & 1-ab)\mqty(-a^{-1} & -b^{-1}\\1 & 1)\\
 &= \mqty((a-b)\bigl(a+a^{-1}\bigr) & 0\\0 & -(a-b)(b+b^{-1}))\\
 &= (a-b)\mqty(a + a^{-1} & 0\\0 & -b-b^{-1}),\\
 U2\tilde{M}_{a_1} U^{*{\mathsf T}} &= \mqty(-a & 1\\-b & 1)\mqty(0 & 1\\-ab & 0)\mqty(-a^{-1} & -b^{-1}\\1 & 1)\\
 &= \mqty(b-a & 0\\0 & a-b) = (a-b)\mqty(-1 & 0\\0 & 1),
 \end{align*}
 so that $UMU^{*{\mathsf T}}$ has the shape that was claimed.
\end{proof}

\begin{Remark}
 Note that in the $q\to1$ limit we have
 \begin{gather*}
 \nabla_k \to \bigl((1-x)\bigl(1-x^{-1}\bigr)\bigr)^{k_1+k_3}\bigl((1+x)\bigl(1+x^{-1}\bigr)\bigr)^{k_2+k_4}\\
\qquad =\bigl(2-x-x^{-1}\bigr)^{k_1+k_3} \bigl(2+x+x^{-1}\bigr)^{k_2+k_4}
 \end{gather*}
 and
 \[
 \tilde{M}\to\mqty(1 & \frac{x+x^{-1}}{2}\\\frac{x+x^{-1}}{2} & 1),\qquad
 U\to\mqty(-1 & 1\\1 & 1)
 \]
 as well as
 \[
 U\tilde{M}U^{*{\mathsf T}}\to \mqty(1-\frac{x+x^{-1}}{2} & 0\\0 & 1+\frac{x+x^{-1}}{2}).
 \]
 This yields the overall matrix weight
 \[
 U\tilde{M}\nabla U^{*{\mathsf T}}\to 2^{k_1+k_2+k_3+k_4} \mqty((1-z)^{k_1+k_3+1}(1+z)^{k_2+k_4} & 0\\0 & (1-z)^{k_1+k_3}(1+z)^{k_2+k_4+1})
 \]
 \big(writing $2z = x + x^{-1}$\big), which corresponds to \cite[equation~(13)]{mvpBC1}.
 Note that by Lemma~\ref{sec:lem-all-affine-root-systems}, we need to set $k_2=k_3=0$
 to emulate an affine root system of type $BC_1$, which leaves the two constants
 $\alpha=k_1$, $\beta=k_4$.
\end{Remark}

Besides the Steinberg basis there happens to be another $A_0$-basis of $A$ that appears
naturally.
\begin{Lemma}\label{sec:lem-km-diagonalise-T1}
 The family $e_1 := 1, e_{s_1} := x^{-1}(1-ax)(1-bx)$ is an $A_0$-basis of $A_J$ ($=A$),
 which diagonalises $T_1$.
\end{Lemma}
\begin{proof}
 We have
 \begin{gather*}
 T_1 \tilde{e}_1 = \tau_1 \tilde{e}_1,\\
 T_1 \tilde{e}_{s_1} = \tau_1 x + (\tau_1-\tau_1^{-1})x^{-1}
 + \tilde{\tau}_1 - \tilde{\tau_1}^{-1}= \tau_1\qty(\frac{1}{ab}x^{-1} + x+x^{-1} - \frac{1}{b} - \frac{1}{a})\\
 \phantom{ T_1 \tilde{e}_{s_1}}{}= \tau_1\qty(\left(x+x^{-1} - \frac{a+b}{ab}\right)\tilde{e}_1 + \frac{1}{ab}\tilde{e}_{s_1}).
 \end{gather*}
 In other words, as a matrix we have
 \[
 T_1 = \tau_1\mqty(1 & x+x^{-1} - \frac{a+b}{ab}\\
 0 & \frac{1}{ab}).
 \]
 This is upper triangular with distinct eigenvalues $1$, $a^{-1}b^{-1}$. Therefore, it can be
 diagonalised with another upper triangular matrix, which transforms $\tilde{e}_1$, $\tilde{e}_{s_1}$ to
$
 1$, $ x^{-1}(1-ax)(1-bx)$.
\end{proof}

\begin{Corollary}
 Let $f,g\in A_0$, then $(fe_1, ge_{s_1})=0$.
\end{Corollary}
\begin{proof}
 Since $T_1$ acts on $A_J$ by means of an $A_0$-endomorphism, its eigenspaces are $A_0$-modules,
 and hence by Lemma~\ref{sec:lem-km-diagonalise-T1}, $fe_1$ is an eigenvector with eigenvalue 1 and $ge_{s_1}$ is an eigenvector
 with eigenvalue $a^{-1}b^{-1}$. Since these are distinct and since $T_1$ is unitary with
 respect to $(\cdot,\cdot)$, the corresponding eigenspaces are orthogonal.
\end{proof}

\begin{Corollary}\label{sec:cor-diagonal-matrix-weight}
 In the expression from Lemma~{\rm\ref{sec:lem-matrix-weight}} of $(\cdot,\cdot)$ in terms of
 a matrix weight for the basis $e_1$, $e_{s_1}$, the matrix weight $M$ is diagonal.
\end{Corollary}

\begin{Lemma}\label{sec:lem-km-coefficients}
 Write \smash{$M=\big(\begin{smallmatrix}d_1 & 0\\0 & d_2\end{smallmatrix}\big)$} for the matrix weight from Corollary~{\rm \ref{sec:cor-diagonal-matrix-weight}}. Then
 \begin{align*}
 d_1 &= \frac{1-ab}{2},\qquad
 d_2 = \frac{1-a^{-1}b^{-1}}{2} \frac{\nabla_{k+l}}{\nabla_k},
 \end{align*}
 where $l=(1,1,0,0)$, so that the labelling $k+l$ corresponds to the coefficients
 $aq$, $bq$, $c$, $d$.
\end{Lemma}
\begin{proof}
 Just like in the proof of Lemma~\ref{sec:lem-weight-matrix-steinberg},
 we have
 \[
 d_1 = \frac{1}{2\delta} (1-s_1) F=\frac{1-ab}{2},\qquad
 d_2 = \frac{1}{2\delta} (1-s_1) (Fhh^*).
 \]
 For the second diagonal entry, we get
 \begin{align*}
 d_2 ={}& \frac{1}{2\delta}(1-s_1) x\bigl(1-ax^{-1}\bigr)\bigl(1-bx^{-1}\bigr)(1-ax)(1-bx)\bigl(1-a^{-1}x^{-1}\bigr)\bigl(1-b^{-1}x^{-1}\bigr)\\
 ={}& \frac{1}{2\delta}\bigl(x\bigl(1-ax^{-1}\bigr)\bigl(1-bx^{-1}\bigr)(1-ax)(1-bx)\bigl(1-a^{-1}x^{-1}\bigr)\bigl(1-b^{-1}x^{-1}\bigr)\\
 & - x^{-1}(1-ax)(1-bx)\bigl(1-ax^{-1}\bigr)\bigl(1-bx^{-1}\bigr)\bigl(1-a^{-1}x\bigr)(1-b^{-1}x)\Big)\\
 ={}&\frac{(1-ax)\bigl(1-ax^{-1}\bigr)(1-bx)\bigl(1-bx^{-1}\bigr)}{2\delta}\Big(x\bigl(1-a^{-1}x^{-1}\bigr)\bigl(1-b^{-1}x^{-1}\bigr)\\
 & - x^{-1}\bigl(1-a^{-1}x\bigr)\bigl(1-b^{-1}x\bigr)\bigr)\\
 ={}& \frac{1-a^{-1}b^{-1}}{2}
 (1-ax)\bigl(1-ax^{-1}\bigr)(1-bx)\bigl(1-bx^{-1}\bigr)\\
 ={}& \frac{1-a^{-1}b^{-1}}{2}
 \frac{\nabla_{k+l}}{\nabla_k}.\tag*{\qed}
 \end{align*} \renewcommand{\qed}{}
\end{proof}

\begin{Corollary}
 The polynomials $(\underline{E_n})_{n\in\Z}$ {\rm(}i.e., the non-symmetric Askey--Wilson polynomials
 expressed in $e_1$, $e_{s_1}${\rm)} form an orthogonal basis with respect to the matrix weight
 \[
 \frac{1}{2}\mqty((1-ab)\nabla_k & 0\\0 & \bigl(1-a^{-1}b^{-1}\bigr) \nabla_{k+l}).
 \]
\end{Corollary}

\begin{Remark}
 In fact, what happened in this subsection is by no means surprising,
 when we consider the contents of \cite[Section~5.8]{macdonald}. Let $\epsilon$
 be the sign character of $W_0$, then according to~\cite[Section~5.8.8]{macdonald}, every $\epsilon$-symmetric polynomial
 lies in $\delta_{\epsilon,k}A_0$ for a particular $\delta_{\epsilon,k}\in A$. It just so happens that $\delta_{\epsilon,k}$
 is a scalar multiple of what we call $e_{s_1}$ here, so that
 we can infer from~\cite[Section~5.8.6]{macdonald} that
$(fe_{s_1}, ge_{s_1})_k \sim (f,g)_{k+l}$
 for $f,g\in A_0$.

 In addition, from \cite[Sections~5.8.10f]{macdonald} (and by comparing leading coefficients), we can conclude that
 \smash{$P_{m,k+l}e_{s_1}=P^{(\epsilon)}_{m+1,k}$}, so that
 writing a polynomial as a two-component vector amounts to determining
 its symmetric and anti-symmetric part (with respect to $T_1$). If we expand the components of $\underline{E_n}$ in terms of
 appropriate symmetric Macdonald polynomials, say as
 \[
 \mqty(\displaystyle \sum_{m\in\N_0} a_{n,m}P_{n,k}\vspace{1mm}\\
 \displaystyle \sum_{m\in\N_0} b_{n,m}P_{n,k+l}),
 \]
 the $Y+Y^{-1}$-eigenvalues dictate that $a_{n,m}=a_n\delta_{\abs{n},m}$
 and $b_{n,m}=b_n\delta_{\abs{n},m+1}$. Consequently, (for~${n\ne0}$), the
 two polynomials $\underline{E_n}$ and $\underline{E_{-n}}$ lie in the space spanned
 by
 \[
 \mqty(P_{n,k}\\0),\qquad \mqty(0\\P_{n-1,k+l}).
 \]
\end{Remark}

\subsection[Example: A\_2 with J=set2]{Example: $\boldsymbol{ A_2}$ with $\boldsymbol{J=\set{2}}$}
This example is the $q$-version of
\cite[Section~7]{mvpVectorValued} (which is linked to matrix spherical
functions of symmetric pairs with restricted root system $A_2$, see
\cite{shimeno}). Here, our affine root system has type~$A_2$, so that $\mathfrak{H}_0$ is generated by $T_1$, $T_2$, and the lattice
$L$ by $\omega_1$, $\omega_2$, the fundamental weights corresponding to $\alpha_1$, $\alpha_2$. There is only one parameter,
$\tau = q(k/2)$.

We consider $J=\set{2}$, leading to the parabolic subgroup $W_J$ generated
by $s_2$ and the subalgebra $\mathfrak{H}_J$ generated by $T_2$. We now
present $A_J$, the $s_2$-invariant polynomials, as a free $A_0$-module.
The construction from~\cite{pittieSteinberg} gives the following basis
$
\tilde{e}_1=1$, $ \tilde{e}_{s_1}=e(\omega_2-\omega_1) + e(-\omega_2)$, $\tilde{e}_{s_2s_1} = e(-\omega_1)$.

Like in the previous subsection, we attempt to get a decomposition by
looking at the action of elements of $\mathfrak{H}_0$ on $A_J$. Those that
let $A_J$ invariant, act as $A_0$-endomorphisms, and therefore can be
expressed as matrices in our existing Steinberg basis.

\begin{Lemma}\label{sec:lem-shi-centraliser-selfadjoint}
 Let $x:= T_1T_2 + T_2T_1 - \bigl(\tau-\tau^{-1}\bigr)(T_1+T_2)$. Then $x\in Z_{\mathfrak{H}_0}(\mathfrak{H}_J)$ and $x^*=x$.
\end{Lemma}
\begin{proof}
 Write $\alpha:=\tau-\tau^{-1}$, then we have
 \begin{align*}
 \begin{split}
& T_2x= T_1 T_2 T_1 + T_1 + \alpha T_2 T_1 - \alpha(T_2T_1
 + 1 + \alpha T_2)= T_1T_2T_1 + T_1 - \alpha^2 T_2 - \alpha,\\
& xT_2= T_1 + \alpha T_1T_2 + T_1T_2T_1
 - \alpha(T_1T_2 + 1 + \alpha T_2)= T_1 + T_1T_2T_1 - \alpha^2 T_2 - \alpha= T_2x,
\end{split}
 \end{align*}
 which shows that $x$ centralises $\mathfrak{H}_J$, the algebra
 generated by $T_2$. For self-adjointness, note that
 \begin{align*}
 x^* &= T_2^{-1}T_1^{-1} + T_1^{-1} T_2^{-1}
 + \alpha(T_1^*+T_2^*)\\
 &= (T_2 - \alpha)(T_1-\alpha) + (T_1-\alpha)(T_2-\alpha)
 + \alpha(T_1 - \alpha + T_2-\alpha)\\
 &= T_2T_1 + T_1 T_2 - 2\alpha T_1 - 2\alpha T_2 + 2\alpha^2
 + \alpha T_1 + \alpha T_2 - 2\alpha^2= x.\tag*{\qed}
 \end{align*} \renewcommand{\qed}{}
\end{proof}

\begin{Lemma}\label{sec:lem-shi-matrix-rep}
 In the basis $\tilde{e}_1$, $\tilde{e}_{s_1}$, $\tilde{e}_{s_2s_1}$, the element
 $x$ from Lemma~{\rm \ref{sec:lem-shi-centraliser-selfadjoint}} has the following matrix
 representation:
 \[
 \mqty(2 & \bigl(\tau^2+1\bigr)m_{\omega_1} & \tau^2m_{\omega_2}\\
 0 & 1-\tau^2 - \tau^{-2} & 0\\
 0 & 0 & 1-\tau^2-\tau^{-2}).
 \]
\end{Lemma}

\begin{Lemma}\label{sec:lem-shi-eigenbasis}
 The family $e_1=1$, $e_{s_1}=h$, $e_{s_2s_1}=h^*$ with
 \[
 h = \bigl(\tau^2+1\bigr)e(\omega_1) - \tau^{-2}(e(\omega_2-\omega_1) + e(-\omega_2))
 \]
 is an $A_0$-basis of $A_J$ that diagonalises $x$.
\end{Lemma}
\begin{proof}
 Note that
 \begin{align*}
& h= \bigl(\tau^2+1\bigr)m_{\omega_1}\tilde{e}_1 - \bigl(\tau^2 + 1 + \tau^{-2}\bigr)\tilde{e}_{s_1},\qquad
 h^* = -\tau^2 m_{\omega_2}\tilde{e}_1 + \bigl(\tau^2 + 1 + \tau^{-2}\bigr)\tilde{e}_{s_2s_1},
 \end{align*}
 so the claimed new basis arises from $\tilde{e}_1$, $\tilde{e}_{s_1}$, $\tilde{e}_{s_2s_1}$ by a triangular matrix whose diagonal entries~$1$
 and $\pm \big(\tau^2+1+\tau^{-2}\big)$ are units of $A_0$. Consequently,
 $e_1$, $e_{s_1}$, $e_{s_2s_1}$ is itself an $A_0$-basis of $A_J$.

 Next, by Lemma~\ref{sec:lem-shi-matrix-rep} $\ker(x-2)$ is the null space of
 \[
 \mqty(0 & \bigl(\tau^2+1\bigr)m_{\omega_1} & \tau^2 m_{\omega_2}\\0 & -1-\tau^2-\tau^{-2} & 0\\0 & 0 & -1-\tau^2-\tau^{-2})
 \]
 and since it is at most 2-dimensional, is spanned by $\tilde{e}_1=e_1$.
 Moreover, the kernel of $x-1+\tau^2+\tau^{-2}$ is the null space of
 \[
 \mqty(\tau^2+1+\tau^{-2} & \bigl(\tau^2+1\bigr)m_{\omega_1} & \tau^2 m_{\omega_2}\\0 & 0 & 0\\0 & 0 & 0),
 \]
 which we can readily see to be spanned by $h$, $h^*$.
\end{proof}

\begin{Corollary}
 Let $f,g\in A_0$, then
$
 (fe_1, ge_{s_1})=(fe_1, ge_{s_2s_1})=0$.
\end{Corollary}
\begin{Corollary}
 By Lemmas~{\rm\ref{sec:lem-shi-eigenbasis}} and~{\rm\ref{sec:lem-shi-centraliser-selfadjoint}}, we have
$
 x fe_1 = fxe_1 = 2fe_1
$
 and
 \[
 x ge_{s_1} = gxe_{s_1} = \bigl(1-\tau^2-\tau^{-2}\bigr)ge_{s_1}
 \]
 and analogously for $e_{s_2s_1}$. We therefore see that $fe_1$, $ge_{s_1}$
 and $fe_1$, $ge_{s_2s_1}$ are pairs of eigenvectors of the self-adjoint
 operator $x$ with distinct eigenvalues. Therefore, they are orthogonal.
\end{Corollary}

We therefore see that our matrix weight from Lemma~\ref{sec:lem-matrix-weight} has the following block shape:
\[
 \mqty(\ast & 0 & 0\\0 & \ast & \ast\\0 & \ast & \ast).
\]
It turns out, there are some relations between these entries that can easily be
determined.

\begin{Lemma}
 We have
 \begin{align*}
 m_{1,1} = \frac{W_0\bigl(\tau^2\bigr)}{6},\qquad
 m_{s_1,s_1} = m_{s_2s_1,s_2s_1},\qquad
 m_{s_1,s_2s_1} = \tau^6 m_{s_2s_1,s_1}^*.
 \end{align*}
\end{Lemma}
\begin{proof} (i) The first claim follows from Lemma~\ref{lem-weight-distributions}\,(iii).
(ii) This follows from the definition of $M$, since $e_{s_1}^*=e_{s_2s_1}$.
(iii) By definition, we have
$
 \Delta_0^{-1} = \tau^3 \bm{c}_{-\alpha_1}\bm{c}_{-\alpha_2}\bm{c}_{-\alpha_1-\alpha_2}$.
 Since $\bm{c}$ is invariant under $*$, we have
 $\Delta_0^{-*} = \tau^{-6} \Delta_0^{-1}$, and hence
 \[
 m_{s_2s_1,s_1}^* = \frac{1}{6}\sum_{w\in W_0}
 w \frac{h^2}{\Delta_0^*}
 = \tau^{-6} \frac{1}{6}\sum_{w\in W_0} w\frac{h^2}{\Delta_0}
 = \tau^{-6} m_{s_1,s_2s_1}.\tag*{\qed}
 \]
\renewcommand{\qed}{}
\end{proof}

This leaves two entries to be determined: $m_{s_1,s_1}$ and $m_{s_1,s_1s_2}$. Our first
observation is that they are both polynomials.

\begin{Lemma}\label{sec:lem-shi-polynomial-weight}
 All $m_{v,v'}$ $\big(v,v'\in W_0^J\big)$ are polynomials.
\end{Lemma}
\begin{proof}
 We proceed similarly to the proof of Lemma~\ref{sec:lem-km-coefficients}: we
 decompose $\Delta_0^{-1}$ in a useful way
 \[
 \Delta_0^{-1} = \tau^3 \prod_{\alpha\in R^+} \bm{c}_{-\alpha,k}
 = \tau^3 \prod_{\alpha\in R^+} \frac{\tau^{-1} e(\alpha/2)-\tau e(-\alpha/2)}{e(\alpha/2)-e(-\alpha/2)}
 =: \frac{F}{\delta},
 \]
 where
 \begin{gather*}
 \delta = \prod_{\alpha\in R^+} (e(\alpha/2)-e(-\alpha/2))= \sum_{w\in W_0} (-1)^{\ell(w)} e(w\rho),\\
F = \tau^3\prod_{\alpha\in R^+} \bigl(\tau^{-1}e(\alpha/2) - \tau e(-\alpha/2)\bigr)= \sum_{w\in W_0} \bigl(-\tau^2\bigr)^{\ell(w)} e(w\rho).
 \end{gather*}
 $\delta$ is the Weyl denominator and every anti-symmetric polynomial (i.e.,
 polynomial that transforms under the sign representation of $W_0$) is divisible by
 it. Consequently, we have
 \[
 m_{v,v'} = \frac{1}{6}\sum_{w\in W_0} w \frac{e_v e_{v'}^* F}{\delta}
 = \frac{1}{6\delta} \sum_{w\in W_0} (-1)^{\ell(w)} w (e_v e_{v'}^* F).
 \]
 Since $e_ve_{v'}^* F\in A$, its anti-symmetrisation is an anti-symmetric polynomial
 and hence divisible by $\delta$. As a consequence, $m_{v,v'}\in A_0$.
\end{proof}

\begin{Lemma}
 The matrix weight $M$ we are considering in this subsection cannot be further reduced.
\end{Lemma}
\begin{proof}
 Assume that it can. By Corollary~\ref{sec:cor-shi-block-structure}, all
 $M_\mu$ can be made diagonal by the same matrix. In particular $M_{2\omega_1}$.
 We will show that this is not possible.

 For that we first find out, which entries of $M_{2\omega_1}$ are definitely zero
 by determining the supports of the coefficients $m_{v,v'}$ ($v,v'\in M_0^J$), i.e.,
 the $W_0$-orbits of $L$ whose coefficient within $m_{v,v'}$ may be nonzero.
 We have the following supports:
 \begin{gather*}
 e_{s_1}e_{s_1}^* = hh^*\colon\ W_0\omega_1 + W_0\omega_2 = W_0\rho\sqcup 0,\\
 e_{s_1}e_{s_1}^*F = Fhh^*\colon\ W_0\rho + (W_0\rho\sqcup 0)
 = W_02\rho\sqcup W_03\omega_1\sqcup W_03\omega_2\sqcup W_0\rho\sqcup 0,\\
 e_{s_1}e_{s_2s_1}^* = h^2 \colon\ W_0\omega_1 + W_0\omega_1 = W_02\omega_1\sqcup W_0\omega_2,\\
 e_{s_1}e_{s_2s_1}^*F = Fh^2\colon\ W_0\omega_1 + (W_02\omega_1\sqcup W_0\omega_2)\\
 \phantom{ e_{s_1}e_{s_2s_1}^*F = Fh^2\colon\ }{}\qquad= W_0(3\omega_1+\omega_2)\sqcup W_0(\omega_1+2\omega_2)\sqcup W_02\omega_1\sqcup W_0\omega_2
 \end{gather*}
 (where we use the polynomial $F$ from the proof of Lemma~\ref{sec:lem-shi-polynomial-weight}).
 Since $m_{s_1,s_1}$ and $m_{s_1,s_2s_1}$ are the anti-symmetrisations of the
 polynomials listed above, divided by $\delta$, we can remove all singular orbits
 since they vanish under anti-symmetrisation, leaving $W_0 2\rho\sqcup W_0\rho$
 in the first case, and $W_0(3\omega_1+\omega_2)\sqcup W_0(\omega_1+2\omega_2)$
 in the second case. Dividing such polynomials by $\delta$ yields the following
 supports
$
 m_{s_1,s_1}\colon W_0\rho\sqcup 0$, $
 m_{s_1,s_2s_1}\colon W_02\omega_1\sqcup W_0\omega_2$.
 This shows that
 \[
M = \mqty(\frac{W_0\bigl(\tau^2\bigr)}{6} & 0 & 0\\0 & a + b m_\rho & \tau^3(cm_{\omega_2} + dm_{2\omega_1})\\0 & \tau^3(c^*m_{\omega_1} + d^*m_{2\omega_2}) & a+bm_\rho)
 \]
 for constants $a,b,c,d\in K$. Therefore,
 \[
 M_{2\omega_1} = d\tau^3 \mqty(0 & 0 & 0\\0 & 0 & 1\\0 & 0 & 0),
 \]
 which can only be brought into diagonal form if $d=0$. To see that $d\ne0$,
 note that of the orbits on which $Fh^2$ is supported, only $W_0(3\omega_1+\omega_2)$ contributes to $d$, i.e., $d\tau^3$ is the coefficient of $e(3\omega_1+\omega_2)$ in the anti-symmetrisation of $Fh^2$. If we expand
$
 h^2 = \sum_\mu a_\mu e(\mu)$,
 we can use the sum expansion of $F$
 \[
 \frac{1}{6}\sum_{w\in W_0} (-1)^{\ell(w)} w\bigl(Fh^2\bigr)
 = \frac{1}{6}\sum_{w,w'\in W_0,\mu} (-1)^{\ell(w)} \bigl(-\tau^2\bigr)^{\ell(w')}
 a_\mu e(w(\mu + w'\rho)).
 \]
 The only way of obtaining $w(\mu + w'\rho)=3\omega_1+\omega_2$ is if
 $w'=w^{-1}$ and $\mu=w'2\omega_1$, thus,
 \[
 d\tau^3 = \frac{1}{6}\sum_{w\in W_0} \tau^{2\ell(w)}
 a_{w^{-1}2\omega_1}
 = \frac{1}{6}\sum_{w\in W_0} \tau^{2\ell(w)} a_{w2\omega_1}
 = \frac{\tau^2+1}{6}\sum_{w\in W_0^I} \tau^{2\ell(w)} a_{w2\omega_1}.
 \]
 Now, the only way of obtaining $2\omega_1$ within $W_0\omega_1+W_0\omega_1$
 is as $\omega_1+\omega_1$, so that
 \[
 a_{2\omega_1} = \bigl(\tau^2+1\bigr)^2,\qquad
 a_{2\omega_2-2\omega_1} = a_{-2\omega_2} = \tau^{-4}.
 \]
 Hence,
 \begin{align*}
 d\tau^3 &= \frac{\tau^2+1}{6}\bigl(\bigl(\tau^2+1\bigr)^2 + \tau^{-4}\bigl(\tau^2+\tau^4\bigr)\bigr)
 = \frac{1+\tau^{-2}}{6} \bigl(\tau^6 + 2\tau^4 + 2\tau^2 + 1\bigr)\\
 &= \bigl(1+\tau^{-2}\bigr)\frac{W_0\bigl(\tau^2\bigr)}{6}\ne0.\tag*{\qed}
 \end{align*} \renewcommand{\qed}{}
\end{proof}

\begin{Corollary}
 If we use the notation from Lemma~{\rm\ref{sec:lem-matrix-weight}}, the polynomials
 \smash{$(\underline{P_{J,\lambda}})_{\lambda\in L_{+,J}}$} are an orthogonal basis of
 $A_0^3$ with respect to the inner product with matrix weight $\nabla M$.
 Here, $M$ is maximally reduced, so we have also found an irreducible $2\times 2$ matrix weight
 on $A_0$.
\end{Corollary}

\subsection*{Acknowledgements}
My thanks to Maarten van Pruijssen for many enlightening and helpful discussions, and to the (anonymous) referees. This research was funded by NWO grant \texttt{OCENW.M20.108}.

\pdfbookmark[1]{References}{ref}
\LastPageEnding

\end{document}